\documentclass{amsart}
\usepackage[top=25mm, bottom=28mm, left=27mm , right=27mm]{geometry}
\usepackage{amsfonts}
\usepackage{amsmath}
\usepackage{amssymb}
\usepackage[applemac]{inputenc}
\usepackage{url}
\usepackage{hyperref}
\usepackage{color}
\usepackage{lipsum}
\usepackage[title]{appendix}
\usepackage{mathtools}
\usepackage{ulem}

\usepackage{derivative}

\usepackage{graphicx}

\def\C{\mathbb{C}}
\def\R{\mathbb{R}}

\def\F{\mathcal{F}}

\def\N{\mathbb{N}}
\def\F{\mathcal {F}}

\def\T{\mathbb{T}}

\def\Heis{{\sf{Heis}} }
\def\Der{{\sf{Der}}}

\def\heis{{\mathfrak{heis}} }

\def\aff{{\mathfrak{aff}} }

\def\Isom{\mathsf{Isom}}

\def\Aff{\sf{Aff}}

\def\Vect{\mathrm{Vect\{}}

\def\Aut{{\sf{Aut}}}

\def\GL{{\sf{GL}}}
\def\O{{\sf{O}}}
\def\SO{{\sf{SO}}}

\def\SL{{\sf{SL}}}

\def\Stab{{\sf{Stab}} }

\def\Lie{{\sf{Lie}} }

\def\ad{\mathsf{ad}}

\def\Ad{\mathsf{Ad}}

\def\Mink{{\sf{Mink}} }

\def\det{{\sf{det}}}

\def\c{{\mathfrak{c}}}

\def\a{{\mathfrak{a}}}
\def\g{{\mathfrak{g}}}

\def\s{{\mathfrak{s}}}

\def\z{{\mathfrak{z}}}
\def\z{{\mathfrak{z}}}

\def\r{{\mathfrak{r}}}

\def\u{{\mathfrak{u}}}

\newtheorem{theorem}{{Theorem}}[section]
\newtheorem{proposition}[theorem]{{Proposition}}
\newtheorem{isom.ext}[theorem]{{Trivial isometric extension}}
\newtheorem{definition}[theorem]{{Definition}}
\newtheorem{lemma}[theorem]{{Lemma}}
\newtheorem{corollary}[theorem]{{Corollary}}
\newtheorem{fact}[theorem]{{\sc Fact}}
\newtheorem{remark}[theorem]{{Remark}}
\newtheorem{question}[theorem]{{Question}}
\newtheorem{observation}[theorem]{{Observation}}
\newtheorem{example}[theorem]{{Example}}
\newtheorem{notations}[theorem]{{Notations}}

\newtheorem{notation}[theorem]{{Notation}}

\definecolor{purple}{rgb}{0.65,0.12,0.94}
\definecolor{forestgreen}{rgb}{0.4,0.64,0.13}

\hypersetup{
	colorlinks=true,
	linkcolor=red,
	filecolor=magenta,      
	urlcolor=blue,
	citecolor=blue
}

\begin{document}

	\title[ ]{On homogeneous plane waves
	}
	
\author [M. Hanounah]{Malek Hanounah}
\address{Institut f\"ur Mathematik und Informatik \hfill\hfill\break\indent
	Walther-Rathenau-Str. 47,
	17489 Greifswald}
\email{malek.hanounah@uni-greifswald.de}

\author [L. Mehidi]{Lilia Mehidi}
\address{Departamento de Geometria y Topologia\hfill\break\indent
	Facultad de Ciencias, Universidad de Granada
	18001 Granada, Spain}
\email{{lilia.mehidi@ugr.es}
	\hfill\break\indent
\url{https://mehidi.pages.math.cnrs.fr/siteweb/}}

\author[A. Zeghib]{Abdelghani Zeghib }
\address{UMPA, CNRS, 
	\'Ecole Normale Sup\'erieure de Lyon\hfill\break\indent
	46, all\'ee d'Italie
	69364 LYON Cedex 07, FRANCE}
\email{abdelghani.zeghib@ens-lyon.fr 
	\hfill\break\indent
	\url{http://www.umpa.ens-lyon.fr/~zeghib/}}

	\date{ \today}
	\maketitle
\begin{center}
\begin{abstract} 
   Plane waves are a special class of Lorentzian spaces with a parallel lightlike vector field.  They are of great importance in Geometry (e.g. Lorentzian holonomy) and in Physics (General Relativity as well as alternative gravity theories). Our contribution in the present paper aims at a rigorous mathematical treatment focusing on completeness of Killing fields, and globality of coordinates.  Equivalence of different approaches to plane waves is by no means easy to handle. We use here cohomogeneity one Heisenberg actions to introduce a point of view from which one can see plane waves as a deformation of Minkowski spacetime. We determine the identity component of the isometry group of a 1-connected non-flat homogeneous plane wave, which establishes a correspondence between these spaces and certain 1-parameter groups of automorphisms of the Heisenberg group.  The extendibility of spacetimes (when incomplete) is a natural, important and delicate question. One of our main results is the proof of the $C^2$-inextendibility of non-flat homogeneous plane waves. We also prove that they are geodesically complete if and only if the lightlike parallel vector field is preserved by the identity component of the isometry group. Finally, we show that a 1-connected homogeneous plane wave admits global Brinkmann coordinates.
\end{abstract}
\end{center}
 
\tableofcontents

\section{Introduction}

 A Lorentzian manifold is called a \textit{Brinkmann manifold} if it admits a parallel lightlike vector field.  They appear in the general context of Lorentzian manifolds with special holonomy, i.e. for which the action of the restricted holonomy group is reducible (having invariant subspaces) and
indecomposable (the metric is degenerate on the invariant subspaces). We denote a Brinkmann manifold by a triple $(\mathbf{X}, g, V)$, where $g$ is the Lorentzian metric, and $V$ the parallel lightlike vector field. 

Since $V$ is parallel, its orthogonal distribution $V^{\perp}$ is invariant by the Levi-Civita connection and is then integrable, hence tangent to a foliation $\mathcal{F}$ having lightlike geodesic leaves (i.e. the induced metric is degenerate). Each leaf of $\mathcal{F}$ inherits an induced connection. A special class of Brinkmann spacetimes are known as  pp-waves: they are defined by the fact that all the $\mathcal{F}$-leaves are flat with respect to the induced connection.

Plane waves are pp-waves which are ``almost symmetric": $\nabla_X R = 0$, for any $X$ tangent to $V^{\perp}$, where $R$ is the Riemannian tensor.  In the so-called local Brinkmann coordinates, the metric takes the form 
\begin{align}\label{Introduction: Eq Brinkmann coordinates}
    2 du dv + S_{ij}(u) x^i x^j du^2 + \sum_{i=1}^{n} (dx^i)^2,
\end{align}
 where $S(u)=(S_{ij}(u))$ is a symmetric matrix. 

The global metrics on $\R^2 \times \R^n$ of the previous form, with $S$  not depending on $u$ and non degenerate (with non-zero determinant), define what we call  \textit{Cahen--Wallach spaces}. They are indecomposable (globally) symmetric plane waves, introduced by Cahen and Wallach \cite{cahen1970lorentzian}.

So pp-waves and plane waves are special classes of Brinkmann spacetimes, obtained with conditions on the curvature. In the opposite direction, (locally) Kundt spacetimes are generalizations of Brinkmann manifolds, defined as Lorentzian manifolds having a codimension one lightlike geodesic foliation.

Plane waves can also be defined by the fact that they admit local Rosen coordinates in which the metric takes the form 
\begin{align}\label{Introduction: Eq Rosen coordinates} 
2 dv du + g_{ij}(u) dx^i dx^j,
\end{align}
where $(g_{ij}(u))$ is a symmetric definite positive matrix. In both Brinkmann and Rosen coordinates, $V=\partial_v$ is the parallel lightlike vector field. 

\subsection{Plane waves as a deformation of Minkowski space}\label{Subsection-Introduction: deformation of Heis action}
 The $(2n+1)$-dimensional Heisenberg group $\Heis_{2n+1}=\R^n \ltimes \R^{n+1}$ is the subgroup of $\mathsf{Aff}(\R^{n+1})$ defined by
$$\Heis_{2n+1}=\left\{\begin{pmatrix}
    1 &\alpha \\
    0 &I_n
\end{pmatrix} | \  \alpha \in \R^{n} \right\}\ltimes \R^{n+1}.$$
This representation of $\Heis_{2n+1}$ in $\mathsf{Aff}(\R^{n+1})$ will be referred to as the usual representation. 
Denote by $A^+\simeq \R^n$ the abelian subgroup of unipotent matrices,  and by  $A^-$ the subgroup $\{0\} \times \R^n$ of the translation part. Let $Z$ be the subgroup of translations $ \R \times \{0\}$ determined by the fixed vector of $A^+$; this subgroup forms the center of the Heisenberg group. Denote their Lie algebras by $\a^+, \a^-$, and $\z$, respectively. 
\begin{notations}
The Lie algebra $\heis_{2n+1}$ of the Heisenberg group equals $\heis_{2n+1}= \mathfrak{a}^+ \oplus \a^- \oplus \mathfrak{z}$ with Lie brackets $$[A,B] = \omega(A,B)z,$$
where  $z$ is a basis for the center $\mathfrak{z}$, and $\omega$ is the standard symplectic form on $\a^+ \oplus \a^- \cong \R^{2n}$. Denote by $X_1, \dots ,X_n$ (resp. $Y_1,\dots,Y_n$) a basis of $\a^+$ (resp. $\a^-$), such that all Lie brackets are zero but $[X_i,Y_i]=z$, for $i=1,\dots,n$.
\end{notations}

\subsubsection{\textbf{The Heisenberg action on the flat Minkowski space}}\label{Introduction_Heis action on Minkowski}
The (flat) Lorentzian space $(\R^{n+2}, g_0= 2du dv + \sum_{i=1}^{n} dx_i^2)$ is called the Minkowski space and denoted by $\Mink^{1,n+1}$.  The group of isometries of $\Mink^{1,n+1}$ is the Poincar\'e group $\O(1, n+1) \ltimes \R^{n+2}$.
Let $\mathsf{SPol}(n)$ be the subgroup of $ \O(1, n+1) \ltimes \R^{n+2}$ preserving a  lightlike constant vector field (this group is introduced in \cite[Section 2]{Article2} and called the special polarized Poincar\'e group).  Up to conjugation by an element in the Lorentz group $\O(1,n+1)$, we can assume that this (constant) lightlike vector field is given by $V=\partial_v$. It has the form $\mathsf{SPol}(n) =  L(\mathsf{SPol}(n)) \ltimes \R^{n+2}$, where the linear part is a semi-direct product $L(\mathsf{SPol}(n)) =\O(n)\ltimes \R^{n}$.  In the basis $(v,x_1,\dots,x_n,u)$ it takes the following form
$$L(\mathsf{SPol}(n)):=\left\{\begin{pmatrix}
    1 & \beta^{\top} & -\frac{\vert \beta \vert^2}{2}\\
0 & A &  - A\beta\\
0 & 0 & 1
\end{pmatrix} |\;\;  A \in \O(n), \beta \in \R^{n}\right\}.$$
The subgroup of $\mathsf{SPol}(n)$ preserving each leaf of the foliation tangent to $V^{\perp}$, given by the level sets $u= cst$, is   $L(\mathsf{SPol}(n))\ltimes \R^{n+1}\cong (\O(n)\ltimes \R^n)\ltimes \R^{n+1}\subset \mathsf{SPol}(n)$, which contains the Heisenberg group $\Heis_{2n+1}=\R^n \ltimes \R^{n+1}$. So this subgroup acts by preserving individually all the lightlike hyperplanes $\R^{n+1} \times \{u\}$, $u \in \R$. We introduce the following groups, which are essential in developing the point of view used here for describing plane waves. 
\begin{definition}\label{Introduction-Definition: Affine unimodular lightlike group}

Let $\R^{n+1}$ with coordinates $(x_0,x_1,\dots,x_n)$  be endowed with the  degenerate quadratic form  $q_0 := x_1^2+\dots+ x_n^2$.

$\mathsf{(1)}$ The group of affine isometries of $q_0$ is given by  
\begin{align*}
	\mathsf{L}(n):&=\left\{\begin{pmatrix}
		\lambda & \alpha^{\top}\\
		0 & A 
	\end{pmatrix} | \ \lambda\in \R^*_+,  \alpha\in \R^{n}, A\in \O(n) \right\}\ltimes \R^{n+1}\\
    &  \cong ((\R \times \O(n)) \ltimes \R^n) \ltimes \R^{n+1}.
\end{align*}
It  will be called the \textbf{affine lightlike group}. The linear subgroup $(\R \times \O(n))\ltimes \R^n$ can be identified with the subgroup of the Lorentz group $\O(1, n+1)$ preserving a lightlike  direction.\\
The affine lightlike group can also be seen as the group  of diffeomorphisms of $\R^{n+1}$ preserving $q_0$ and the (natural) flat affine connection of the affine space. 

$\mathsf{(2)}$ The subgroup $\mathsf{L}_\u(n):= (\O(n) \ltimes \R^n) \ltimes \R^{n+1}$ of $\mathsf{L}(n)$ preserving the lightlike constant vector field $\partial_{x_0}$ will be called the \textbf{affine unimodular lightlike group}.\medskip

A manifold modeled on $(\mathsf{L}(n),\R^{n+1})$  $($resp. $(\mathsf{L}_\u(n),\R^{n+1}))$ will be said to have an affine  lightlike geometry (resp. affine unimodular lightlike geometry). \bigskip
\end{definition}

Using the vocabulary in Definition \ref{Introduction-Definition: Affine unimodular lightlike group}, we see that $\Heis_{2n+1}$ acts on each hyperplane $\R^{n+1} \times \{u\}$ of the Minkowski space by preserving the induced affine unimodular lightlike geometry.  So for every $u \in \R$, this action defines a faithful representation
		$$\theta_u: \Heis_{2n+1} \to \mathsf{L}_\u(n):= \Aut(\R^{n+1}, \nabla, \frac{\partial}{\partial v}, dx_1^2+\dots+dx_n^2),$$
where $\nabla$ denotes the affine flat connection of $\R^{n+1}$. 

Let us take a closer look at this action. At the level $u=0$, $\Heis_{2n+1}=\R^n \ltimes \R^{n+1}$ acts through the usual action, where $\R^{n+1}$ acts by translations and $\lambda \in \R^{n}$ acts linearly by a unipotent transformation 
$\begin{pmatrix}
    1 &\lambda \\
    0 & I_n
\end{pmatrix}$.
At the level $u$, $\R^{n+1} \subset \Heis_{2n+1}$ acts by translations and $\lambda \in \R^{n+1}$ acts with an affine transformation $A(\lambda)+T_u(\lambda)$ depending on $u$, where
$$A(\lambda)=\begin{pmatrix}
    1 &\lambda \\
    0 & I_n
\end{pmatrix} \mathrm{\;\;and\;\;} T_u(\lambda)=\begin{pmatrix}
    -\frac{\vert \lambda\vert^2}{2} u \\
    -\lambda u
\end{pmatrix}.$$ 
\noindent So the $\R^{n+1}$-action by translations is constant with respect to $u$, but the $\R^n$-action depends on $u$.  
In fact, one can write the $\theta_u$-action of $\mathsf{Heis}_{2n+1}$ on $\R^{n+1} \times \{u\}$ as a deformation of its $\theta_0$-action on $\R^{n+1} \times \{0\}$  by means of an automorphism of $\Heis_{2n+1}$ that depends on $u$. Namely, $$\theta_u =  \theta_0 \circ  \mathfrak{P}_u^{\Mink},$$ 
where $\mathfrak{P}_u^{\Mink} \in \Aut(\Heis_{2n+1})$ belongs to the automorphism group of $\Heis_{2n+1}$ and is given by
$$d_0 \mathfrak{P}_u^{\Mink}= \exp(uL), \
L=\begin{pmatrix}
   0 & 0 & 0\\
   I_n & 0 & 0\\
   0 & 0 & 0
\end{pmatrix},
$$ where the derivation $L \in \Der(\heis_{2n+1})$ is written in the decomposition $\heis_{2n+1}=\a^+ \oplus \a^- \oplus \z$. 
	And the global action of $\mathsf{Heis}_{2n+1}$ on the Minkowski space can be written as follows:
	\begin{align*}
		a_{\mathfrak{P}^\mathsf{Mink}}:  \mathsf{Heis}_{2n+1}\times \mathbb{R}^{n+1} \times \mathbb{R} &\to  \mathbb{R}^{n+1} \times \mathbb{R} \\
		(h, (v,x,u) ) & \mapsto (\theta_0 \circ \mathfrak{P}_u^{\mathsf{Mink}}(h) (v,x), u).
	\end{align*} 

The deforming curve $\mathfrak{P}_u^{\Mink}$ is a $1$-parameter group in $\Aut(\Heis_{2n+1})$. 

\subsubsection{\textbf{Plane waves as curves in $\Aut(\Heis_{2n+1})$}}\label{Section: Plane waves as curves in Aut(Heis)} 
Plane waves can be thought of as a generalization as well as a deformation of the Minkowski space. As explained above, the Minkowski space $\Mink^{1,n+1}$ admits an isometric action of $\Heis_{2n+1}$, acting transitively on the leaves of some codimension $1$ foliation by (parallel) affine lightlike hyperplanes $\R^{n+1} \times \{u\}$, $u \in \R$. It acts on a fixed hyperplane $\R^{n+1} \times \{u\} \cong \R^{n+1}$ by an affine transformation that depends on $u$ (in particular, this is not a product action). This determines a curve of affine maps of $\R^{n+1}$ depending on $u$, that can be naturally parameterized by some curve in $\Aut(\Heis_{2n+1})$.  We will use this point of view of cohomogeneity one $\Heis_{2n+1}$-actions through parametrizing curves in $\Aut(\Heis_{2n+1})$ to give another description of plane waves, which make them appear as a deformation of the Minkowski space. This will be precised in the following paragraph. Let us however note that the computations leading to the equivalence between this point of view and the original definition are not trivial.

\subsection*{Plane wave case}
The leaves of the $\F$-foliation of a plane wave (and more generally, of a pp-wave) are flat and lightlike, with a parallel lightlike vector field $V$ tangent to them. Hence, they inherit a unimodular affine lightlike geometry. As a result, the Lie algebra of Killing fields tangent to the $\F$-foliation and commuting with $V$ is (faithfully) represented in $\mathfrak{o}(n) \ltimes \heis_{2n+1}$, the Lie algebra of $\mathsf{L}_{\u}(n)$. What is the image of this representation? It turns out that, in the case of plane waves, the image contains the full $\heis_{2n+1}$ algebra, and this property characterizes plane waves among pp-waves (see Section \ref{Section: Synthetic description of plane waves}). This is a consequence of the well known fact (at least in the indecomposable case) that plane waves admit an isometric infinitesimal action of the Heisenberg algebra, whose action preserves individually the leaves of $\mathcal{F}$ (see for instance \cite{ehrlich1992gravitational} in dimension $4$, and \cite{Blau} in the indecomposable case). 
This infinitesimal action of the Heisenberg algebra does not necessarily extend to a global action of the Heisenberg group (see the discussion in the next paragraph).

We will see in Section \ref{Section: Synthetic description of plane waves} that similarly to the Minkowski case, for which the $\Heis_{2n+1}$-action determines some $1$-parameter group $\mathfrak{P}^{\Mink} \in \Aut(\Heis_{2n+1})$,  a plane wave determines  certain curves in $\Aut(\Heis_{2n+1})$ that encode the $\Heis_{2n+1}$ local action. We formulate the following question, which we answer in Section \ref{Section: Synthetic description of plane waves}:
\begin{question}
    \label{Intro: Problem Heis action}
Let $I \subset \R$ be an interval, and consider $\mathbf{Y}=\R^{n+1} \times I$. Define an action of $\Heis_{2n+1}$ on $\mathbf{Y}$ as follows
		\begin{align*}
			a_\mathfrak{P} :  \Heis_{2n+1}\times \R^{n+1} \times I &\to  \R^{n+1} \times I\\
			(h, (y, u) ) & \mapsto (\theta \circ \mathfrak{P}_{u}(h) (y), u),
		\end{align*} 
		where  
		\begin{itemize}
			\item $\mathfrak{P}: I \to \Aut(\Heis_{2n+1})$ is some curve,
			\item and $\theta: \Heis_{2n+1} \to \mathsf{L}_{\u}(n)$ is some faithful representation into the affine unimodular lightlike group. 
		\end{itemize}
Assume the $a_\mathfrak{P}$-action preserves a Lorentzian metric on $\mathbf{Y}$, which Lorentzian metrics do we obtain?
\end{question}

In the question above, we allow $\mathfrak{P}$ to be any curve. This $a_{\mathfrak{P}}$-action of $\Heis_{2n+1}$ is a deformation of the  $a_{\mathfrak{P}^\Mink}$-action obtained in the case of Minkowski space.
As we will see, not every such action preserves a Lorentzian metric on $\mathbf{Y}$. However, when it does, the metric turns out to be a plane wave. A natural question then arises: how can we characterize all such curves that arise from plane wave metrics on $\mathbf{Y}$? This characterization is provided in Section \ref{Section: Synthetic description of plane waves}. Moreover,  $\Heis_{2n+1}$ acts locally isometrically on a plane wave via the restriction of such an $a_{\mathfrak{P}}$-action.

\subsection*{Plane waves in standard form} When Brinkmann coordinates exist globally on a plane wave $\mathbf{X}:=\R^{n+1} \times I$, for an open interval $I \subset \R$ parameterizing the space of $\F$-leaves, we refer to it as a plane wave \textbf{in standard form} (this definition appears for instance in \cite[Section 1]{Leis-Schlie}). In Section \ref{Section: Synthetic description of plane waves}, we will see that a plane wave in standard form admits a global isometric $a_{\mathfrak{P}}$-action of $\mathsf{Heis}_{2n+1}$ as above.  We write explicitly the $\Heis_{2n+1}$-action in the (global) Brinkmann coordinates in the Appendix.  In particular, we will be able to describe a plane wave in standard form  as a single curve in $\Aut(\Heis_{2n+1})$. 

\subsection{Isometry groups of $1$-connected homogeneous plane waves}

Recall that the Lie algebra of Killing fields of an indecomposable plane wave contains the Heisenberg algebra, whose action is locally transitive on the $\F$-leaves.  Actually, it turns out that the result remains true in the decomposable case (this is considered in Section \ref{Section: Synthetic description of plane waves}). Thus, general plane waves have already local cohomogeneity $1$.

In the homogeneous case, Blau and O'Loughlin \cite{Blau} determine the full isometry algebra of an indecomposable  homogeneous plane wave.
They start with a plane wave in local (Brinkmann or Rosen) coordinates, and ask when it is (infinitesimally) homogeneous. Through a direct analysis of the Killing equations in these coordinates, they identify the isometry algebra as being generated by Killing fields satisfying specific bracket relations, and prove that it contains a Heisenberg algebra. It appears then that the Killing algebra contains an $\R$-extension of the Heisenberg algebra, though the action of the one-parameter group of (infinitesimal) automorphisms on the Heisenberg algebra is not visible.
To explicitly describe $1$-connected homogeneous plane waves as Lorentzian homogeneous spaces, which is the aim of Theorem \ref{Introduction-Theorem: isometry group} below, one has to determine their  isometry groups. This problem in general requires understanding which local isometries extend globally. 
This is achieved here by analyzing the globalization of local isometries generated by the Heisenberg algebra, and deducing the full isometry group relying on two key points: first, that the foliation $\F$ carries a unimodular affine lightlike structure (see Definition \ref{Introduction-Definition: Affine unimodular lightlike group}) along its leaves, and second, that it also admits an affine transverse structure (Section \ref{Section: Affine structure on the space of leaves}), both of which are preserved by the isometries preserving $\F$.

Studying whether the action of $\heis_{2n+1}$ integrates into an action of the Heisenberg group (i.e., whether the local isometries extend globally) amounts to studying the completeness of Killing fields tangent to the $\F$-foliation, which in turn involves  the (geodesic) completeness of the $\F$-leaves of such manifolds.  
This is investigated in Section \ref{Section: completeness of F-leaves}.
For this, we use the fact that the leaves of the $\F$-foliation of a pp-wave inherit a unimodular affine lightlike geometry.
It is a well known fact that if $M$ is a $(G,X)$-manifold such that the model space $X$ is complete for some $G$-invariant connection, the completeness of $M$ as a $(G,X)$-manifold is equivalent to its geodesic completeness for the connection it inherits from $X$. In the homogeneous case, we prove the following: 

\begin{proposition}\label{Introduction-Proposition: homog implies F-complete} 
1) The $\F$-leaves of a non-flat homogeneous plane wave $(\mathbf{X},g,V)$ of dimension $n+2$ are $(\mathsf{L}_\u(n), \R^{n+1})$-complete, hence geodesically complete.\\
2) A homogeneous pp-wave (possibly flat) such that $V$ is complete has complete $\F$-leaves.
\end{proposition}

As a consequence, it follows that the infinitesimal action of the Heisenberg algebra on a $1$-connected homogeneous plane wave integrates into a Lie group isometric action of $\Heis_{2n+1}$. Moreover, as in the flat case, $\Heis_{2n+1}$ acts on $\F$ by preserving each individual leaf of $\F$, and the action can be parameterized by a $1$-parameter group in $\Aut(\Heis_{2n+1})$. With this in hand, we are able to determine the identity component of the isometry groups of $1$-connected homogeneous plane waves (without using the Killing equation). It turns out that these groups contain an $\R$-extension of the Heisenberg group. The one-parameter groups of automorphisms on the Heisenberg group are described explicitly; they fall into two families: one for which the plane wave is complete, and one for which it is incomplete.
Finally, we show that if $(\mathbf{X},g,V)$ is a $1$-connected non-flat homogeneous plane wave and $G$ is the identity component of its isometry group, then any $G$-invariant Lorentzian metric on $X$ is isometric to $g$.
This is stated in the following theorem

\begin{theorem}\label{Introduction-Theorem: isometry group}
Let $(\mathbf{X},g,V)$ be a $1$-connected non-flat homogeneous plane wave. Then
\begin{enumerate}
\item \textbf{Isometry group.} The connected component of $\Isom(\mathbf{X})$ (which is of finite index in $\Isom(\mathbf{X})$) is isomorphic to $G_\rho=(\R \times K) \ltimes \Heis_{2n+1}$,  with the $\R$-action given by $\rho(t)=e^{t L}$, where $L \in \Der(\heis)$ is such that 
     \begin{itemize}
        \item[(a)] \
          $L=\begin{pmatrix}
           F & B & 0\\
           I & F & 0\\
           0 & 0 & 0
        \end{pmatrix}$,  or  
        \item[(b)] \ 
          $L=\begin{pmatrix}
            I+F & B & 0\\
            I & F & 0\\
            0 & 0 & 1
           \end{pmatrix}$
     \end{itemize}
    with respect to the decomposition  $\heis_{2n+1}=\a^+\oplus \a^- \oplus \z$, where $F$ is antisymmetric and $B$ is symmetric. Moreover, $K$ is the connected (closed) subgroup of $\O(n)$ with Lie algebra $\mathfrak{k}=\{E \in \mathfrak{o}(n), [E,F]=[E,B]=0 \}$, acting on $\Heis_{2n+1}$ trivially on the center and by the standard action on $A^+$ and $A^-$. \medskip
    
  Finally, $\mathbf{X}$ identifies with $X_\rho=G_\rho/K\ltimes A^+$. Moreover, the codimension $1$ foliation $\mathcal{F}$ tangent to $V^\perp$ is given by the (left) action of $K \ltimes \Heis_{2n+1}$, and is preserved by the action of $G_\rho$.\medskip

\item \textbf{Homogeneous plane waves vs $1$-parameter groups in $\Aut(\Heis_{2n+1})$.} Let $G_\rho$ be as in item (1), and consider the homogeneous space $X_\rho=G_\rho/I$, with $I=K \ltimes A^+$. Then $X_\rho$ admits a unique (up to isometry of $X_\rho$) $G_\rho$-invariant Lorentzian metric, which is necessarily a plane wave metric (which may be flat).   \medskip 
\item \textbf{Completeness.} If the $\R$-action is as in case (a), the plane wave is geodesically complete. Otherwise, i.e. case (b), they are incomplete. 
\end{enumerate}
\end{theorem}
\begin{remark}[The flat case]\label{Remark: isometry group of flat PW}
As shown in \cite{duncan1989homogeneous}, a $1$-connected flat homogeneous Lorentzian space is isometric to either the whole Minkowski space or a half-Minkowski space bounded by a lightlike hyperplane $P$, i.e., with a degenerate induced metric. Both spaces are homogeneous flat plane waves, with $S=0$ in (\ref{Introduction: Eq Brinkmann coordinates}). The isometry group in the first case is the Poincar\'e group. In the second case, it is given by the subgroup of the Poincar\'e group that preserves the boundary lightlike hyperplane $P$ and the connected components of $\R^{n+2} \smallsetminus P$, or, equivalently, it is the subgroup preserving a lightlike direction tangent to $P^\perp$. This group is nothing but the affine lightlike group $\mathsf{L}(n)$ introduced in Definition \ref{Introduction-Definition: Affine unimodular lightlike group}; it has the form of the groups obtained in Theorem \ref{Introduction-Theorem: isometry group} above.
\end{remark}
As already mentioned in Subsection \ref{Subsection-Introduction: deformation of Heis action}, a similar description as in Item \textbf{(2)} of Theorem \ref{Introduction-Theorem: isometry group} can be done with non-homogeneous plane waves in standard form, but this time, the action of $\Heis_{2n+1}$ is parameterized by a curve in $\Aut(\Heis_{2n+1})$ (which is a $1$-parameter group if and only if the space is homogeneous).  Let us observe that the assumption of being in standard form ensures that the $\heis_{2n+1}$-action integrates into an action of the Lie group. Otherwise, one can ask the following questions for general $1$-connected plane waves 
\begin{question}[\textbf{Global $\Heis_{2n+1}$-action vs $\F$-completeness}]\label{Intro-Question 1}
Is there an equivalence between the two following properties? 
   \begin{enumerate}
		\item[(a)] The infinitesimal action of the Heisenberg algebra integrates to a global $\Heis_{2n+1}$-action. 
		\item[(b)] The leaves of $\mathcal{F}$ are complete. 
	\end{enumerate}  
\end{question}
\begin{question}\label{Intro-Question 2}
 One then also asks whether $1$-connected plane waves with complete $\F$-leaves are exactly the plane waves in standard form.    
\end{question}
We do not have an answer to the above questions yet. However, establishing a positive answer to Question \ref{Intro-Question 1} would allow to extend the picture in item (2) to the family of $1$-connected plane waves with complete $\F$-leaves.
\subsection{Existence of global Brinkmann coordinates} In \cite{Blau}, assuming the existence of global Brinkmann coordinates, Blau and O'Loughlin classified homogeneous plane waves in standard form. However, it is not clear whether any $1$-connected homogeneous plane wave necessarily admits such coordinates.  Theorem \ref{Intro: Theorem Brinkmann} below fills this gap.
Indeed, as an application of Theorem \ref{Introduction-Theorem: isometry group}, we prove the existence of global Brinkmann coordinates on a $1$-connected homogeneous plane wave (Section \ref{Section: Global Brinkmann coordinates}), providing an affirmative answer to Question \ref{Intro-Question 2} in the homogeneous case. The spaces we obtain coincide with those found in \cite{Blau}.

\begin{theorem}\label{Intro: Theorem Brinkmann}
   Let $\mathbf{X}$ be a $1$-connected  homogeneous plane wave. Then $\mathbf{X}$ has global Brinkmann coordinates.
\end{theorem}
\subsection{$C^2$-inextendibility}
The question of the extendibility of a spacetime is an important and difficult question, which arose from the resolution of Einstein equations in general relativity. Choquet-Bruhat and Geroch \cite{choquetbruhat} proved that if a  globally hyperbolic spacetime $M$ is a solution to   the Cauchy problem for  Einstein equations, then there is a (unique) maximal extension of $M$ as a solution of the Einstein equations, which is globally hyperbolic.
  On the other hand,   Penrose's  \textit{strong cosmic censorship conjecture} states that ``generically", this extension is also maximal as a Lorentzian spacetime.  
Plane waves are   not in general globally hyperbolic, and they do not care natural foliations by Cauchy hypersurfaces. Instead, they are     foliated by lightlike hypersurfaces (the $\F$-foliation).  
It turns out that in dimension $3$, all plane waves are conformally flat, and that their conformal development in the Einstein universe is a non-trivial conformal extension. The question of conformal extendibility is considered in \cite{Ross}, where they mention that some plane waves admit conformal extensions and others may not. Here, we are interested in  the metric extendibility of plane waves.
It turns out that the lightlike foliation of a plane wave has geodesic leaves, so, the extendibility question can be raised first along the foliation $\F$.  One then asks
\begin{question}
Let $\mathbf{X}$ be a $1$-connected plane wave. Can we always embed $\mathbf{X}$ in a plane wave which is $\F$-complete?  $\F$-maximal ?
\end{question}
In the homogeneous (non-flat) case, we show that the $\F$-leaves are complete.  Moreover, as mentioned in Theorem \ref{Introduction-Theorem: isometry group}, a $1$-connected homogeneous plane wave for which the parallel lightlike vector field is preserved by the (identity component of) isometry group is complete. Otherwise, it is incomplete. In the latter case, we prove that it is $C^2$-inextendible.
\begin{theorem}\label{Introduction-Theorem: maximality}
Let $\mathbf{X}$ be a $1$-connected non-flat homogeneous plane wave, then $\mathbf{X}$ is $C^2$-inextendible, i.e. there is no $C^2$-embedding of such a plane wave as a proper open subset of some Lorentzian manifold.  
\end{theorem}

\subsection*{Organization of the paper}
In Section \ref{Section: Synthetic description of plane waves}, we write the $\Heis_{2n+1}$-action explicitly in Brinkmann and Rosen coordinates, following a computation performed in the Appendix, and we describe plane waves as curves in $\Aut(\Heis_{2n+1})$.
In Section \ref{Section: Affine structure on the space of leaves}, we see that the space of $\F$-leaves of a $1$-connected Brinkmann manifold has an affine structure, inducing a representation of the identity component of the isometry group in $\mathsf{Aff}(\R)$. In Section \ref{Section: completeness of F-leaves}, when the space is $1$-connected and homogeneous, we study the kernel of this representation, consisting of all the isometries acting trivially on the space of $\F$-leaves. We prove Proposition \ref{Introduction-Proposition: homog implies F-complete}, and obtain that $\Heis_{2n+1}$ is contained in the kernel. Section \ref{Section: Isometry group} is devoted to give the full form of the identity component of the isometry group (Theorem \ref{Introduction-Theorem: isometry group}). As an application, we show in Section \ref{Section: Global Brinkmann coordinates} that such spaces admit global Brinkmann coordinates (Theorem \ref{Intro: Theorem Brinkmann}). Finally, Section \ref{Section: C^2 maximality} is devoted to the proof of Theorem \ref{Introduction-Theorem: maximality} on $C^2$-inextendibility.

 \subsection*{Acknowledgement}
We would like to thank Ines Kath for helpful discussions and comments on this paper, and for providing us with the global Brinkmann coordinates (Section \ref{Section: Global Brinkmann coordinates}). We also thank the referee for the thorough evaluation and valuable comments and suggestions, that helped improve the quality of the presentation.
\\
The second author is partially supported by the grants, PID2020-116126GB-I00\\
(MCIN/ AEI/10.13039/501100011033), and the framework IMAG/ Maria de Maeztu,\\ 
CEX2020-001105-MCIN/ AEI/ 10.13039/501100011033.

\section{Synthetic description of plane waves}\label{Section: Synthetic description of plane waves}
In the first part of the section, we explicit the $\Heis_{2n+1}$ local action on a plane wave, in both Rosen and Brinkmann coordinates.  In both cases, when the coordinates exist globally on the plane wave,  there is a global isometric $\Heis_{2n+1}$-action on it, and the metric determines some curve in $\Aut(\Heis_{2n+1})$ that encodes this action; we will write this curve explicitly in both cases. 
The second part of this section is devoted to Question \ref{Intro: Problem Heis action}. 

\subsection{Local $\Heis_{2n+1}$-actions}
\subsubsection{\textbf{In Rosen coordinates}}
In this subsection,  $\theta_0: \Heis_{2n+1} \to \mathsf{Aff}(\R^{n+1})$ is the  action of $\Heis_{2n+1}=\R^n \ltimes \R^{n+1}$ on $\R^{n+1}$, in which $\R^{n+1}$ acts by translation and $\lambda \in \R^{n}$ acts linearly by a unipotent transformation
$ \begin{pmatrix}
		1 & \lambda  \\
		0 & I_n   \\
\end{pmatrix}.$\\ 
Consider $\mathbf{Y}=\R^{n+1} \times I$ with coordinates  $(v,x=(x_1,\dots,x_n),u)$, equipped with the following metric
	\begin{equation} \label{1}
		g = 2 dv du + g_{ij}(u) dx^i dx^j.
	\end{equation}  
The vector field $V=\partial_v$ is parallel, and $g$ is the metric of a plane wave in Rosen coordinates. Suppose $0 \in I$, and denote by $Q(u)$ the symmetric positive definite matrix $(g_{ij}(u))$.\\
	
	\noindent The Heisenberg group $\Heis_{2n+1}= \R^n \ltimes \R^{n+1}$ acts isometrically on $(\mathbf{Y},g)$ as follows: 
	
	$\bullet$ $\R^{n+1}$ acts by translation on the $(v,x)$-coordinates, and trivially on the $u$-coordinate.  
	
	$\bullet$ $\R^n$ acts linearly by unipotent transformations, namely, $\lambda=(\lambda_1,\dots,\lambda_n) {}^\top \in \R^n$ acts as follows	
	\begin{align*}
		\lambda \cdot  
		\begin{pmatrix}
			v \\
			x\\
			u
		\end{pmatrix}&= 
		\begin{pmatrix}
			\begin{pmatrix}
				1 & \lambda {}^\top   \\
				0 & I_n   \\
			\end{pmatrix}
			\begin{pmatrix}
				v \\
				x
			\end{pmatrix}+
			\begin{pmatrix}
				\int_{0}^{u} \alpha \\
				\int_{0}^{u} \beta 
			\end{pmatrix}\\
			u
		\end{pmatrix},
	\end{align*}
 where $\alpha$, $\beta$ are given by $$\beta(u,\lambda):= -Q(u)^{-1}( \lambda) \mathrm{\;\;and\;\;} \alpha(u, \lambda):=- \frac{1}{2}  \lambda {}^\top Q(u)^{-1} \lambda.$$  
So the $\R^n$-factor acts trivially on the $u$-coordinate, and the action on each subspace $\R^{n+1} \times \{u\}$ is affine, with the translation part depending on $u$.\\
\\
Observe that on the level $u=0$, the $\Heis_{2n+1}$-action on $\R^{n+1}\times \{0\}$ coincides with the $\theta_0$ action given in the beginning of this paragraph. On another $u$-level, 
the $\Heis_{2n+1}$-action on $\R^{n+1} \times \{u\}$ can be seen as a deformation of the $\theta_0$ action by applying an automorphism of $\Heis_{2n+1}$. More precisely, define a curve $\mathfrak{P}: I \to \Aut(\Heis_{2n+1})$ in the following way
	\begin{align*}
		&\mathfrak{P}_{u_{\vert A^- \times Z}} = Id		\\
		&\mathfrak{P}_u(\lambda) = (\lambda, \int_{0}^{u} \alpha, \int_{0}^{u} \beta) \in A^+ \times Z \times A^-, \mathrm{\;\;for\;\;} \lambda \in A^+	
	\end{align*}
	For every $u \in I$, the map $\mathfrak{P}_u$ extends uniquely to an automorphism of $\Heis_{2n+1}$. Then the action of $\Heis_{2n+1}$ on $\mathbf{Y}=\R^{n+1} \times I$ can be written as
	\begin{align*}
		a_{\mathfrak{P}} :  \Heis_{2n+1}\times \R^{n+1} \times I &\to  \R^{n+1} \times I\\
		(h, (y, u) ) & \mapsto (\theta_0 \circ \mathfrak{P}_{u}(h) (y), u). 
	\end{align*}
 
\subsubsection{\textbf{In Brinkmann coordinates}}
In this subsection,  we use another common realization of $\Heis_{2n+1}=\C^n \times \R$ as a central extension of the abelian group $A=\R^{2n}$ by $Z=\R$ defined by
$$(a,z)\cdot (a',z')=(a+a', z+z'+\frac{1}{2}\omega(a,a')),$$
for $z,z' \in Z$ and $a,a' \in A$, and where $\omega$ is the standard symplectic form on $\R^{2n}$. This induces another representation $\theta_1: \Heis_{2n+1} \to \mathsf{Aff}(\R^{n+1})$ for which the center $Z=\R$ acts on $\R^{n+1}$ by translation along the $v$-coordinate, and $(\lambda', \lambda) \in \R^{2n}$ acts through the unipotent affine transformation 	
	$ \left(\begin{pmatrix}
		1 & \lambda'   \\
		0 & I_n   \\
	\end{pmatrix},
	\begin{pmatrix}
		\frac{ \langle\lambda, \lambda' \rangle}{2}   \\
		\lambda   \\
	\end{pmatrix}\right) 
	\in \SL_{n+1} (\R) \ltimes \R^{n+1}$. Define the subspaces $A^+:= \R^n \times \{0\}, A^-:= \{0\} \times \R^n$ of $A$. \\
\\
Consider $\mathbf{Y}=\R^{n+1} \times I$ with coordinates $(v,x=(x_1,\dots,x_n),u)$, equipped with the metric
\begin{equation} \label{Eq1: Brinkmann coordinates}	
	g = 2 dv du + x {}^\top S(u) x\,  du^2 + \Sigma_{i=1}^n dx_i^2,
\end{equation}  
where $S(u)$ is a symmetric matrix. This is the metric of a plane wave in Brinkmann coordinates, and $V:=\partial_v$ is a parallel lightlike vector field.\\

Let $S: u \in I \mapsto S(u)$ be the curve of symmetric matrices defining the metric $g$. We denote by $C_{\O(n)}(S)$ the subgroup of $\O(n)$ consisting of elements that commute with $S(u)$ for all $u \in I$. 

\begin{proposition}\label{Proposition: computing G_0 in Brinkmann coordinates} 
Let $(\mathbf{Y},g,V)$ be a plane wave in standard form (\ref{Eq1: Brinkmann coordinates}), and $G_0$ the subgroup of $\Isom(\mathbf{Y},g)$ acting trivially on the $u$-coordinate. Then $G_0$ is isomorphic to $C_{\O(n)}(S) \ltimes \Heis_{2n+1}$, where  $C_{\O(n)}(S)$ acts on $A^-$ and $A^+$ by its standard action, and trivially on the center $Z$ of $\Heis_{2n+1}$. Moreover, the vector field $V$ is a generator of the center $\z$ of $\heis_{2n+1}$. 
\end{proposition}
	\noindent The action of an element $\varphi \in G_0=C_{\O(n)}(S) \ltimes \Heis_{2n+1}$ is computed in the Appendix. It has the following form
	
	\begin{align}\label{Action of G_0 in Brinkmann coordinates}
		\varphi  
		\begin{pmatrix}
			v \\
			x\\
			u
		\end{pmatrix}&=
		\begin{pmatrix}
			v - \langle\alpha'(u), Ax + \frac{1}{2} \alpha(u)\rangle + c  \\
			Ax+\alpha(u)  \\
			u
		\end{pmatrix}\\
		&=
		\begin{pmatrix}
			\begin{pmatrix}
				1 & - \alpha'(u) {}^\top   \\
				0 & I_n   \\
			\end{pmatrix}
			\begin{pmatrix}
				1 & 0   \\
				0 & A   \\
			\end{pmatrix}
			\begin{pmatrix}
				v \\
				x
			\end{pmatrix}+
			\begin{pmatrix}
				-\frac{1}{2} \langle \alpha'(u), \alpha(u) \rangle \\
				\alpha(u) 
			\end{pmatrix} +	
			\begin{pmatrix}
				c \\
				0
			\end{pmatrix}\\
			u
		\end{pmatrix}
	\end{align}
	where $ c \in \R$, $\alpha$ is a solution of the differential equation
	\begin{align}\label{Equation alpha}
		\alpha''(u) = S(u) \alpha(u), \;\forall u \in I,
	\end{align}
	and $A \in \O(n)$ is a constant orthogonal matrix that commutes with $S(u)$ for every $u \in I$.\medskip
	
Denote by $E$ the $2n$-dimensional vector space of solutions of $(\ref{Equation alpha})$. Let $\{e_1,\dots,e_n\}$ be the canonical basis of $\R^n$.  Let $u_0 \in I$, and denote by $\alpha_i$ (resp. $\beta_i$) the solution such that $$\alpha_i(u_0)=0 \ \text{and} \ \alpha'_i(u_0) = e_i \ \text{(resp. $\beta_i(u_0)=e_i$ and $ \beta'_i(u_0) = 0$)}.$$  This identifies $E$ to $\R^{2n}$ by sending $\alpha \in E$ to $(\alpha'(u_0), \alpha(u_0)) \in \R^{2n}$.\medskip

\begin{observation}\label{Observation: rho_1 action of H}
Let $h:=(\lambda', \lambda, c) \in \Heis_{2n+1}$, and fix $u_0 \in I$. Let $\varphi \in G_0$ be the map for which $A=I_n$ and $\alpha$ is the solution of $(\ref{Equation alpha})$ satisfying $\alpha'(u_0)=\lambda', \alpha(u_0)=\lambda$. Then, the action of $\varphi$ on the level $\R^{n+1} \times \{u_0\}$ is nothing but the $\theta_1$-action of $h$ on $\R^{n+1}$.
\end{observation} 
\begin{proof}[Proof of Proposition \ref{Proposition: computing G_0 in Brinkmann coordinates}]
Denote by $H$ the subgroup of elements of $G_0$ with $A=I$. We will prove that $H$ is isomorphic to the Heisenberg group of dimension $2n+1$. First, observe that since $G_0$ acts trivially on the space of $u$-leaves, the action of an element in it is completely determined by its action on some $u$-level. Fix some level $u_0$, and define $\Bar{H}:=\{\varphi_{\vert \{u=u_0\}}, \varphi \in H \}$. So $H$ is isomorphic to $\Bar{H}$. Now define a correspondence from $\Heis_{2n+1}=\C^n \times \R$ to $\Bar{H}$ that sends an element $h:=(-\lambda',\lambda,c) \in \Heis_{2n+1}$ to the affine map $\varphi_{u_0}: \R^{n+1} \to \R^{n+1}$,
$$ \varphi_{u_0}(v,x)= \begin{pmatrix}
			\begin{pmatrix}
				1 & - \lambda'^\top   \\
				0 & I_n   \\
			\end{pmatrix}
			\begin{pmatrix}
				v \\
				x
			\end{pmatrix}+
			\begin{pmatrix}
				-\frac{1}{2} \langle \lambda', \lambda \rangle \\
				\lambda 
			\end{pmatrix} +	
			\begin{pmatrix}
				c \\
				0
			\end{pmatrix}
		\end{pmatrix}.
$$
Clearly, $\varphi_{u_0}$ is the restriction to the $u_0$-level of an element in $H$, namely the map $\varphi$ for which $\alpha$ is the the solution of $(\ref{Equation alpha})$ satisfying $\alpha'(u_0)=\lambda', \alpha(u_0)=\lambda$. This correspondence is one-to-one, and is clearly a group morphism by Observation \ref{Observation: rho_1 action of H}. Now, the structure of $G_0$ can be seen by composing two elements of the form (\ref{Action of G_0 in Brinkmann coordinates}). And the fact that $V$ generates the center $\z$ of $\heis_{2n+1}$ follows from Observation \ref{Observation: rho_1 action of H}.
\end{proof}
Here again, the $\Heis_{2n+1}$-action on the level $\R^{n+1} \times \{u\}$ can be obtained by deforming a $\theta_1$-action on some level $\R^{n+1} \times \{u_0\}$, $u_0 \in I$, by applying an automorphism of $\Heis_{2n+1}$ depending on $u$. Namely, define a curve $\mathfrak{P}: I \to \Aut(\Heis_{2n+1})$ by 
\begin{align*}
		 &\mathfrak{P}_{u_{\vert Z}} = Id \\
		 &\mathfrak{P}_{u}(\lambda',\lambda,0) = (-\beta'(u), \beta(u),  0 ) \in A^+ \times A^- \times Z, \mathrm{\;\; for \;\;} (\lambda', \lambda) \in A^+ \times A^-,
\end{align*}
	where $\beta$ is the solution of Equation $(\ref{Equation alpha})$ satisfying $\beta(u_0) =\lambda, \beta'(u_0) = \lambda'$.  For every $u \in I$, the map $\mathfrak{P}_u$ extends uniquely to an automorphism of $\Heis_{2n+1}$. Then	the action of $\Heis_{2n+1}$ on $\mathbf{Y}=\R^{n+1} \times I$ can be written as
	\begin{align*}
		a_\mathfrak{P} :  \Heis_{2n+1}\times \R^{n+1} \times I &\to  \R^{n+1} \times I\\
		(h, (y, u) ) & \mapsto (\theta_1 \circ \mathfrak{P}_{u}(h) (y), u). 
	\end{align*}

\begin{remark}
		$\mathfrak{P}$ is a (local) $1$-parameter group in $\Aut(\Heis_{2n+1})$ if and only if the differential equation (\ref{Equation alpha}) is autonomous, i.e. $u \mapsto S(u)$ is constant, which amounts to saying that $(\mathbf{Y},g,V)$ is a Cahen-Wallach space. 
\end{remark}
	
\subsection{Plane waves as curves in $\Aut(\Heis_{2n+1})$}
Any plane wave admits local Rosen and Brinkmann coordinates. Thus, by the previous paragraphs, given a plane wave $\mathbf{X}$ of dimension $n+2$,  for any point $p \in \mathbf{X}$, there exists a neighborhood $U$ of $p$ with an embedding $i: U \to \mathbf{Y} :=\mathbb{R}^{n+1} \times I$ such that the local action of $\mathsf{Heis}_{2n+1}$ on $U$ corresponds to the restriction to $i(U)$ of an $a_{\mathfrak{P}}$-action of $\mathsf{Heis}_{2n+1}$ on  $\mathbf{Y}$, as described above. The $a_{\mathfrak{P}}$-actions given above are deformations of the  $a_{\mathfrak{P}^\mathsf{Mink}}$-action on the Minkowski space described in Paragraph \ref{Introduction_Heis action on Minkowski}.

We  consider here more general  $a_{\mathfrak{P}}$-deformations of  the $a_{\mathfrak{P}^\mathsf{Mink}}$-action. Let  $\mathbf{Y} :=\mathbb{R}^{n+1} \times I$, with $0 \in I$. Assume  that $\mathsf{Heis}_{2n+1}$ acts on $\mathbf{Y}$ such that, on a level $\mathbb{R}^{n+1} \times \{u\}$, for $u \in I$, this action is given by a faithful representation $$\theta_u: \mathsf{Heis}_{2n+1} \to   \mathsf{L}_\mathfrak{u}(n) = \mathsf{O}(n) \ltimes \Heis_{2n+1}$$ 
of $\mathsf{Heis}_{2n+1}$ into the lightlike unimodular group $ \mathsf{L}_\mathfrak{u}(n)$ (see Definition \ref{Introduction-Definition: Affine unimodular lightlike group}), depending on $u$.

\begin{lemma}\label{Lemma: rho(Heis)=Heis}
  Let $\varphi: H \to \mathsf{L}_\u(n)=\O(n)\ltimes \Heis_{2n+1}$ be a morphism with a discrete kernel, with $H$ a nilpotent Lie group of dimension $2n+1$. Then $\varphi(H)=\Heis_{2n+1}$ (in particular $\varphi$ is faithful).
\end{lemma}
\begin{proof}
  Denote by $p$ the projection from $\O(n)\ltimes \Heis_{2n+1}$ to $\O(n)$. Let $k:= \dim p(\varphi(H))$. Since $\varphi(H)$ is nilpotent, and the dimension of a maximal abelian subgroup of $\O(n)$ is $\lfloor \frac{n}{2}\rfloor$, we have  $k\leq \lfloor \frac{n}{2}\rfloor$. Computing the brackets $[(A,h), h']$, where $(A,h)\in \varphi(H)$ and $h'\in \varphi(H)\cap \Heis_{2n+1}$, together with the fact that $\varphi(H)$ is nilpotent, shows that the action of $p(\varphi(H))$ on $\varphi(H)\cap \Heis_{2n+1}$ is trivial. It follows that $p(\varphi(H))$ embeds in $\O(E^\perp)$, where $E:= (\varphi(H)\cap \Heis_{2n+1})^{\mathsf{o}}$. But $\dim E= 2n+1-k$; as a consequence, $k\leq \lfloor \frac{k}{2}\rfloor$, which is a contradiction unless $k=0$. 
\end{proof}

By Lemma \ref{Lemma: rho(Heis)=Heis}, for every $u \in I$,  we have $\theta_u(\Heis_{2n+1})=\Heis_{2n+1}$. It follows that $$\theta_u= \theta_0 \circ \mathfrak{P}_u$$ for some $\mathfrak{P}_u \in \Aut(\Heis_{2n+1})$. This defines a  curve  $\mathfrak{P}: u \in I \mapsto   \Aut(\Heis_{2n+1})$, such that the action of $\Heis_{2n+1}$ on $\mathbf{Y}$  is given by 
\begin{equation}\label{Deformed action}
	    \begin{aligned}
		a_\mathfrak{P} :  \Heis_{2n+1} \times \R^{n+1} \times I &\to  \R^{n+1} \times I\\
                (h, (y, u) ) & \mapsto (\mathfrak{P}_u(h) \cdot y, u) 
	\end{aligned}
	\end{equation}
where $\cdot$ is the  $\theta_0$-action on $\R^{n+1} \times \{0\}$. This motivates the following question:
\begin{question}[Question \ref{Intro: Problem Heis action}]
 Fix a (faithful) representation $\theta_0: \Heis_{2n+1} \to \mathsf{L}_\u(n) =\O(n) \ltimes \Heis_{2n+1}$. Let $I \subset \R$ be an interval such that $0 \in I$.
Define a curve $\mathfrak{P}: u \in I \mapsto   \Aut(\Heis_{2n+1})$. And consider $a_{\mathfrak{P}}$  the $\Heis_{2n+1}$ action on $\mathbf{Y}:= \R^{n+1} \times I$ given by (\ref{Deformed action}). For which curves $\mathfrak{P}$ the    $a_{\mathfrak{P}}$-action preserves a plane wave metric on $\mathbf{Y}$? 
\end{question}
The rest of this paragraph is devoted to answering the above question. 
Denote the coordinates $((v,x),u) \in \R^{n+1} \times I$. 

\begin{proposition}\label{Proposition: action preserves Lorentz metric => plane wave}
A Lorentzian metric on $\mathbf{Y}$ preserved by $a_{\mathfrak{P}}$ is a plane wave metric. Moreover, a vector field  $V$ generating the action of the center $\z$ of $\heis_{2n+1}$ is lightlike and parallel, and the $\F$-leaves tangent to $V^\perp$ are given by the $u$-levels. 
\end{proposition}
\begin{proof}
	Let $p_u=(0,u) \in \mathbf{Y}$. 
For the $\theta_u$-action of $\Heis_{2n+1}$ on $\R^{n+1} \times \{u\} \subset \mathbf{Y}$, the isotropy group of $p_u$ is $Is(p_u)=\theta^{-1}_u(A^+)$, where $$A^+=\left\{\begin{pmatrix}
1 & \lambda^{\top}\\
0 & I_n 
\end{pmatrix} | \ \lambda\in \R^n  \right\}.$$ 
It acts on  $\R^{n+1} \times \{u\}$ through its representation $\theta_u(Is(p_u))=A^+$. 
Denote by $e_0 \in T_{p_u} \R^{n+1}$ the fixed vector of $A^+$. A quadratic form on $\R^{n+1} \times \{u\}$ invariant under the $\theta_u$-action of $\Heis_{2n+1}$ corresponds to a scalar product $q$ on $T_p \R^{n+1}$ invariant under the action of the isotropy representation. It is easily seen that such $q$ satisfies $q(e_0,e_0)=0$ and $e_0^\perp=\R^{n+1}$. 
Denote by $V$ the vector field generating the center $\z$ of $\heis_{2n+1}$ such that $V(p)=e_0$. Since the $\z$-action commutes with the $\heis_{2n+1}$-action, the vector field $V$ is everywhere lightlike on $\R^{n+1} \times \{u\}$. 
Assume now that the $a_{\mathfrak{P}}$ action preserves a Lorentzian metric on $\mathbf{Y}$. It follows from the observations above that the induced metric on $\R^{n+1} \times \{u\}$ is degenerate for every $u \in I$, and that $V$ is everywhere lightlike on $\mathbf{Y}$. 

Now, fix $u \in I$. The subgroup  $\theta^{-1}_u(A^- \times Z)$ of the Heisenberg group is abelian, and acts freely in some neighborhood of any point in the level $\R^{n+1} \times \{u\}$, generating  $n+1$ commuting (local) Killing vector fields that span the lightlike distribution in this neighborhood. It follows, using \cite[Theorem 3]{Leis}, that the metric in this neighborhood is that of a plane-wave, and the lightlike vector field $V$ is parallel. 
\end{proof}
\begin{notation}
 Denote by $X_1,\dots,X_n$ (resp. $Y_1,\dots,Y_n$) the canonical basis of $\a^+$ (resp. $\a^-$), where $X_i:=
 \left(\begin{pmatrix}
     0 & e_i^\top\\
     0 & 0_n
 \end{pmatrix}, 0\right)
 $ and $Y_i= \left(0_{n+1}, e_i \right)$ $(0_k$ is the zero matrix of size $k \times k)$. All Lie brackets are zero but $[X_i,Y_i]=z$, for $i=1,\dots,n$.   
\end{notation} 
The next proposition characterizes $a_{\mathfrak{P}}$-actions preserving a Lorentzian metric on $\mathbf{Y}$. 
Fix $p_u:=(0,u) \in \mathbf{Y}$, $u \in I$. The isotropy at $p_u$ is given by $Is(p_u)=\mathfrak{P}_u^{-1} \circ \theta_0^{-1}(A^+)$.
Let $\lambda \in A^+$, then $\mathfrak{P}_u^{-1} \circ \theta_0^{-1}(\lambda) \in Is(p_u)$ acts on $\R^{n+1} \times \{u+\varepsilon\}$, $\varepsilon \geq 0$,  via an affine map
	\begin{equation}\label{Eq: action of A+ through affine map}
 \theta_0 \circ P_{\varepsilon} \circ \theta_0^{-1}(\lambda) =	\left(\begin{pmatrix}
		1 & a_{\varepsilon}(\lambda)  \\
		0 & I_n 
	\end{pmatrix},
 \begin{pmatrix}
		c_{\varepsilon}(\lambda)  \\
		b_{\varepsilon}(\lambda) 
	\end{pmatrix}\right) \in \Heis_{2n+1},
\end{equation}
where $P_{\varepsilon}:=\mathfrak{P}_{u+\varepsilon} \circ \mathfrak{P}_{u}^{-1} \in \Aut(\Heis_{2n+1})$. Of course this map also depends on $u$, but since $u$ is fixed, we indicate only the dependence on $\varepsilon$ in the matrices, for simplicity.\\ 
\\
Note that for any $\varepsilon \geq 0$, $b_{\varepsilon}(\lambda)$ is linear in $\lambda$, so it determines a $1$-parameter family of linear maps $L: \varepsilon \geq 0 \mapsto L(\varepsilon)=b_{\varepsilon} \in \mathrm{Hom}( A^+, A^-)$. Moreover, using that $P_\varepsilon \in \Aut(\Heis_{2n+1})$ for all $\varepsilon \geq 0$, one shows easily that the matrix $L'(0)$ represented in the canonical basis $(X_1,\dots,X_n)$ and $(Y_1,\dots,Y_n)$ of $\mathfrak{a}^+$ and $\mathfrak{a}^-$ respectively, is symmetric (this can also be seen directly by Lemma \ref{Lemma: derivations of heis}). Set $D:= L'(0)$. Finally, we denote by $C$ the vector $C :=c'(0)$, where $c: \varepsilon \geq 0 \mapsto (c_{\varepsilon}(X_1),\dots, c_{\varepsilon}(X_n))$. 

\begin{proposition}\label{Proposition: Existence of a_theta invariant Lorentz metric}
	$a_{\mathfrak{P}}$ preserves a Lorentzian metric $g$ on  $T \mathbf{Y}_{|\R^{n+1} \times \{u\}}$  if and only if $D$ is definite. Moreover, $g=f(u)g_0$, with $f(u):= g(\partial_u,V)$,  and $g_0$ is determined by $D$ up to choosing the scalar product $\langle \partial_u,\partial_u \rangle_{(0,u)}$, and satisfies $g_0(\partial_u,V)=1$. 
\end{proposition}	
\begin{proof}
Since $\Heis_{2n+1}$ has transitive action on the level $\R^{n+1} \times \{u\}$, $a_{\mathfrak{P}}$ preserves a Lorentzian metric $g$ on  $T \mathbf{Y}_{|\R^{n+1} \times \{u\}}$ if and only if the isotropy at $p_u=(0,u)$ preserves some Lorentzian quadratic form $q$ at the tangent space $T_{p_u} \mathbf{Y}$. Let $\theta^{-1}_u(\lambda) \in Is(p_u)$, with $\lambda \in A^+$. Write its action as
	\begin{align*}
		\phi  
		\begin{pmatrix}
			v \\
			x\\
			u+\varepsilon
		\end{pmatrix}&= 
		\begin{pmatrix}
			\begin{pmatrix}
				1 & a_\varepsilon(\lambda)   \\
				0 & I_n   \\
			\end{pmatrix}
			\begin{pmatrix}
				v \\
				x
			\end{pmatrix}+
			\begin{pmatrix}
				c_\varepsilon(\lambda) \\
				b_\varepsilon(\lambda)
			\end{pmatrix}\\
			u+\varepsilon
		\end{pmatrix}.
	\end{align*}
Then
 	\begin{align}\label{Eq-proof: isotropy representation}
		d_{p_u} \phi = 
		\begin{pmatrix}
				1 & \lambda & c(\lambda)   \\
				0 & I_n & b(\lambda)  \\
                0 & 0   & 1
		\end{pmatrix},
	\end{align}   
 where $c(\lambda) := \frac{d}{d\varepsilon}_{\vert_{\varepsilon=0}} c_\varepsilon(\lambda)$ and $b(\lambda) := \frac{d}{d \varepsilon}_{\vert_{\varepsilon=0}}  b_\varepsilon(\lambda)$. 
Let $q$ be a quadratic form at $T_{p_u} \mathbf{Y}$ preserved by $d_{p_u} \phi$. Set $e_0 = \partial_v, e_1 = \partial_{1},\dots, e_n = \partial_{n}, e_{n+1} = \partial_u$ at $p_u$.
Then necessarily
 \begin{itemize}
     \item $q(e_0,e_0)=0$
     \item $q(e_0,e_i)=0$ for any $i=1,\dots,n$
     \item $H(b(\lambda)) = - q(e_0,e_{n+1}) \lambda$, where $H:=(q(e_i,e_j))_{i,j}$, which gives $$H=-q(e_0,e_{n+1}) D^{-1}$$
     \item $D^\top A= -q(e_0,e_{n+1})B$, where $A=(q(e_i,e_{n+1}))_i$ ($i=1,\dots,n$) and $B:=C - \frac{1}{2}   (D_{ii})_i$,.
 \end{itemize}
We see that $q= q(e_0,e_{n+1}) q_0$, where $q_0$ is determined by $D$, up to choosing $q(e_{n+1},e_{n+1})$. It is Lorentzian if and only if $D$ is positive (resp. negative) definite and $q(e_0,e_{n+1}) <0$ (resp. $>0$). 

To finish the proof, we will show that $g$ factorizes through  $g(\partial_u, \partial_v)$, and that the latter only depends on $u$. First, observe that by Proposition \ref{Proposition: action preserves Lorentz metric => plane wave}, the vector field $\partial_v$ is everywhere lightlike for $g$, and $\partial_v^\perp=\Vect \partial_v,\partial_1,\dots,\partial_n\}$. 
Now, the scalar product at another point $(x,v,u)$ of $\R^{n+1} \times \{u\}$ is the push-forward of the scalar product $q$ at $T_{p_u} \mathbf{Y}$ by means of an element in $\theta^{-1}_u(A^- \times Z)$. This element acts on $\R^{n+1} \times \{u\}$ via some affine map as in (\ref{Eq: action of A+ through affine map}), whose differential has then the same form as (\ref{Eq-proof: isotropy representation}). It then sends $\partial_v$ to $\partial_v$, and $\partial_u$ to $\partial_u + X$, with $X \in \Vect \partial_v,\partial_1,\dots,\partial_n\}$. So the scalar product $g(\partial_v, \partial_u)$ is constant on a fixed $u$-level, as stated. 
\end{proof}

\begin{remark}
We proved that a Lorentzian metric preserved by the $a_{\mathfrak{P}}$-action is a plane wave. One can ask whether plane waves obtained this way admit global Brinkmann coordinates, i.e. coincide with the plane waves in standard form. As mentioned in Questions \ref{Intro-Question 1} and \ref{Intro-Question 2} of the introduction, one may ask even a more general question: whether existence of a global isometric action of $\Heis_{2n+1}$ on a $1$-connected plane wave implies existence of global Brinkmann coordinates. 
\end{remark}

\section{Affine structure on the space of leaves}\label{Section: Affine structure on the space of leaves}

\subsection{Affine structure on the space of leaves of a Brinkmann manifold}\label{Subsection: Affine structure on the space of leaves}
Let $(\mathbf{X},g,V)$ be a $1$-connected Brinkmann spacetime. Let $\mathcal{F}$ be the foliation determined by $V^{\perp}$. 
Let $G = \Isom^\mathsf{o}(\mathbf{X}, g,\R V)$ be the identity component of the subgroup of the isometry group of $\mathbf{X}$ preserving the line field $\R V$ (since $G$ is connected, it actually preserves the direction field defined by $V$). This means that for any $f \in G$, the push-forward $f_*$ sends $V$ to $a V$, for some positive function $a$. Since $V$ is parallel,  $a$ must be a constant in  $ \R^*_+$. 
	
The action of $G$ is uniquely determined by its infinitesimal action 
\begin{align*}
		\sigma : \mathfrak{g} &\to \Gamma(T\mathbf{X})\\
		Y &\mapsto Y_\mathbf{X}, \mathrm{\;\;with\;\;} Y_\mathbf{X}(p) =  \left. \frac{d}{dt} \right|   _{t=0} \exp(tY) \cdot p    
\end{align*}
where $\mathfrak{g}=\Lie(G)$ identifies then with the Lie algebra of complete Killing  fields of $\mathbf{X}$ preserving the line field $\R V$,  i.e., those whose flows send a leaf of  $V$ to a leaf of $V$.
	
\begin{proposition}\label{Proposition: affine structure on ksi}
	 Let $\mathbf{X}$ be a $1$-connected Brinkmann spacetime.	The space of leaves $\xi:=\mathbf{X}/\F$ is a $1$-dimensional  (non necessarily Hausdorff) manifold, with an affine structure preserved by any element of $G$.
\end{proposition}
\begin{proof}
	The foliation $\mathcal{F}$  is defined by the non-singular closed $1$-form $\omega=g(V,.)$, so is clearly transversely affine  (we can refer to  \cite{blumenthal}, which defines the more general transversally homogeneous foliations). 
	Let $p \in \mathbf{X}$, one can write $\omega = du$, where $u: \mathbf{X}\to \R, u(p) =0,$ is a submersion defining the leaves of $\mathcal{F}$.  The local sections of $u$ give local charts on $\xi$, and define a $1$-dimensional manifold structure on it. 
	The map $u$  induces a local diffeomorphism from $\xi$ to an interval $I$ of $\R$ (which is a global diffeomorphism when $\xi$ is Hausdorff), and defines a global affine parameter on $\xi$. 
    
    Now, for any $f \in G$, $u \circ f : \mathbf{X}\to \R$ is another submersion defining another affine parameter on $\xi$ that can be deduced from $u$ by an affine transformation of $\R$. Indeed, since $f$ preserves $\R V$,  there exists $a \in \R^*_+$ such that $df(V) = a^{-1} V$. This yields $f^*\omega = a \omega$, hence $d(u \circ f) = a du$. So for any $f \in G$, there exists $(a,b) \in \R^* \ltimes \R$ such that $u \circ f = a u +b$, where $b=u(f(p))$   (in particular, the linear part of the induced affine transformation is determined by the action of $f$ on $V$). This means that the transverse affine structure on $\F$ is preserved by $G$.  
\end{proof}
\begin{notation}\label{Notation: submersion u}
    In all this section, $u$ denotes the submersion (defined up to translation) $u \in C^{\infty}(\mathbf{X},\R)$ such that $du=\omega$, defining the transverse affine structure of $\F$ (it is introduced in the proof of Proposition \ref{Proposition: affine structure on ksi} above). 
\end{notation}
	
	It follows that $G$ acts as an automorphism of the affine space $\xi$. So we have a representation 
	$$\pi: G \to \mathsf{Aff}(\xi) =\mathsf{Aff}_I(\R),$$
    where $ \mathsf{Aff}_I(\R)$ is the subgroup of $\mathsf{Aff}(\R)$ preserving $I$. 	And we have an exact sequence:
\begin{equation}\label{Exact sequence on G}
	1 \rightarrow  G_0 \rightarrow G \rightarrow  \pi(G) < \mathsf{Aff}(\xi) \rightarrow 1 
\end{equation}
where $G_0 = \ker \pi$ is a normal subgroup of $G$ that acts trivially on the space of leaves $\xi$, with Lie algebra 
\begin{equation}\label{g_0}
  \mathfrak{g}_0 = \{ Y \in \mathfrak{g}, g(Y_\mathbf{X},V)=0 \}.  
\end{equation}
\bigskip
The real line has three non-isomorphic affine structures: 
	\begin{enumerate}
		\item  The real line $\R$ with automorphism group $\R^* \ltimes \R$, 
		\item  The half line $]0, + \infty [$ with automorphism group the subgroup of homotheties $\R_+^*$,
		\item $]-1,1[$ with automorphism group $\{\pm \mathrm{Id}\}$.
	\end{enumerate}
	In the first two cases, the automorphism group acts transitively on the manifold, but not in the third case. We deduce the following fact:
\begin{fact}
	If $\mathbf{X}$ is $G$-homogeneous, then the space of $\mathcal{F}$-leaves $\xi$ is Hausdorff, and its affine structure is either of type $(1)$ or $(2)$. Moreover:
\begin{itemize}
		\item	 In case (1), the subgroup of $G$ consisting of isometries preserving $V$ is mapped surjectively by $\pi$ onto the translation subgroup of $\Aff(\mathbb{R})$. Furthermore, it contains $G_0$, the kernel of $\pi$, which consists exactly of the isometries preserving $V$ and fixing a leaf of $\mathcal{F}$.
		\item  In case (2), the kernel $G_0$ coincides with the subgroup of isometries preserving $V$.
	\end{itemize}
\end{fact}

\subsection{The action of Killing fields}\label{Subsection 3.2}
Let $\mathbf{X}$ be a $1$-connected locally homogeneous Brinkmann manifold.  
A classical result \cite{nomizu} says that any locally defined Killing vector field extends coherently to $\mathbf{X}$. 
Let $\Bar{\mathfrak{g}}$ be the Lie algebra of Killing fields of $\mathbf{X}$ preserving the line field $\R V$. The Killing fields need not be complete, hence the use of a notation different from that of the previous paragraph. Given a Killing field $K$, $K \in \Bar{\g}$  is equivalent to $[K,V]$ being collinear to $V$, or, since $V$ is parallel, to $[K,V]=\alpha V$, for some constant $\alpha \in \R$. \medskip

 The flow of $V$ acts isometrically on $\mathbf{X}$: denote by $\mathfrak{z}$ the Lie subalgebra generating the action of $V$.  
\begin{fact}
	$\mathfrak{z}$ is an ideal of $\Bar{\mathfrak{g}}$.
\end{fact}
The local flow of a Killing field in $\Bar{\mathfrak{g}}$ preserves the affine structure of $\xi$, and induces therefore  a (local) $1$-parameter group of affine transformations of $I$, hence an element in $\aff(\R)$. So we have a representation 
\begin{align}\label{Eq: representation of Lie algebra in aff(R)}
\overline{\pi}: \Bar{\mathfrak{g}} \to \aff(\R).   
\end{align}
To write this explicitely, let $\{X,T\}$ be the standard basis of $\aff(\R)$, with $[X,T] = T$. Let $e^{t X}$ (resp. $e^{t T}$) be the $1$-parameter group  of homotheties $u \mapsto e^{t} u$ (resp. translations $u \mapsto u + t$) of $\R$.  Recall the submersion $u \in C^{\infty}(\mathbf{X},\R)$ such that $du=\omega$, defining the transverse affine structure of $\F$ (see Notation \ref{Notation: submersion u}). 
The invariance of the transverse affine structure by $K \in \Bar{\mathfrak{g}}$ yields the existence of $(\alpha,\beta) \in \R^* \times \R$ such that $K \cdot u = \alpha u+\beta$, and $\overline{\pi}$ maps $K$ to $\alpha X + \beta T \in \aff(\R)$. Indeed, since $\omega$ is closed, using Cartan's formula, one can see that  $[K,V] = \alpha V$ for some $\alpha \in \R$  if and only if $\mathcal{L}_K \omega = \alpha \omega$, if and only if $K \cdot u = \alpha u+\beta$, for some $\beta \in \R$, which yields the element $\alpha X + \beta T \in \aff(\R)$ associated to $K$.  The latter also reads $g(K,V) = \alpha u+\beta$. 
\vspace{0.2cm}
\subsection{The homogeneous case}
Let $\mathbf{X}$ be a $1$-connected homogeneous Brinkmann manifold, admitting a unique parallel lightlike vector field $V$ up to scale (this is the case for instance when $\mathbf{X}$ is indecomposable). Then $\mathbf{X}$ is $G$-homogeneous, with $G=\Isom^{\mathsf{o}}(\mathbf{X},g,\R V)$. 
\begin{notations}
 Throughout the paper, the Lie subgroup of $G$ tangent to $\mathfrak{z}$ will be denoted by $Z$. It is a closed subgroup of $G$  contained in the center of $G_0$.  This is consistent with the notation $Z$ of the center of the Heisenberg group, since when $X$ is a plane wave, $Z$ corresponds exactly to the center of the Heisenberg group (Proposition \ref{Proposition: computing G_0 in Brinkmann coordinates}).
\end{notations}

\begin{proposition}\label{Propositiion: Isometry group of Homogeneous Brinkman}
		Let $\mathbf{X}$ be a $1$-connected homogeneous Brinkmann manifold, admitting a unique parallel lightlike vector field $V$ (up to scale). Then
		\begin{itemize}
			\item[$\mathsf{(1)}$] If $\dim \pi(G) = 1$, then  $ G= \R \ltimes_{\rho} G_0$, for a morphism $\rho: \R \to \Aut(G_0)$. In this case, the affine space $\xi$ is complete exactly when $Z$ is in the center of $G$.  
			\item[$\mathsf{(2)}$] If $\dim \pi(G) = 2$, then $G=\R \ltimes (\R \ltimes G_0)$.  Let $G_1$ be the normal subgroup of elements preserving $V$. If the action of $G_0$ on the leaves of $\mathcal{F}$ is not transitive, the $G_1$-orbits in $\mathbf{X}$ define a codimension $1$ foliation of $\mathbf{X}$ by Brinkmann submanifolds. 
		\end{itemize}	
\end{proposition}	
\begin{proof}
	 $\mathsf{(1)}$  If $\mathbf{X}$ is homogeneous, $\xi$ is Hausdorff, affinely equivalent to either $\R$ or $\R^*_+$. 
	 Suppose $\dim \pi(G)=1$. If $\xi = \R$,  then $\pi$ represents $G$ in $\mathsf{Aff}(\R)$, and	 $\pi(G)$ is a Lie subgroup of $\mathsf{Aff}(\R)$ that acts transitively on $\R$. Such a Lie subgroup must contain the $1$-parameter group of translations (generated by $T \in \aff(\R)$). Let $K \in \mathfrak{g}$ such that $d_e \pi(K) = T$. Then $K_\mathbf{X}$ is a complete Killing field acting transitively on $\xi$. Its flow is a $1$-parameter subgroup of $G$ that defines a splitting of the exact sequence   $1 \to G_0 \to G \to \R \to 1$. We get $G=\R \ltimes G_0$ for some $\rho: \R \to \Aut(G_0)$. Furthermore, $[K,V]=0$, so $Z$ is in the center of $G$. 
	 Now, if $\xi = \R^*_+$, then $\pi: G \to \R^*_+$ is surjective, i.e. $\pi(G)$ is the group of homotheties $\R_+^*$. As before, we get $G$ as an extension of $G_0$ by a $1$-parameter group of $\Aut(G_0)$. Here, $Z$ is not in the center of $G$.\\
    $\mathsf{(2)}$ Now suppose $\dim \pi(G)=2$. Since $G_1$ is the inverse image by $\pi$ of the normal subgroup of translations of $\mathsf{Aff}(\R)$, it is normal in $G$, and we have an exact sequence $1 \to G_1 \to G \to \R \to 1$. Hence $G= \R \ltimes G_1$. On the other hand, $G_0$ is a normal subgroup of $G_1$ and we have $G_1=\R \ltimes G_0$. Finally $G=\R \ltimes (\R \ltimes G_0)$ as stated. 
\end{proof}
Let $K \in \mathfrak{g}$ be the Killing field on $\mathbf{X}$, transverse to $V^{\perp}$, generated by the action of the $\mathbb{R}$ factor in $G$. Let $\phi^t$ be the global flow of $K$. We have a global diffeomorphism $\psi: \R \times \mathcal{N}_0 \to \mathbf{X}$, $\psi(t,p) = \phi^t(p)$, where $\mathcal{N}_0$ is a leaf of $\mathcal{F}$. The $G$-action on $\mathbf{X}$ reads as a $\R \ltimes_{\rho} G_0$ action on $\R \times \mathcal{N}_0$ as follows: 
\begin{align}\label{Eq: R ltimes G_0 action on RxN_0}
(\R \ltimes_{\rho} G_0) \times \R \times \mathcal{N}_0  &\to \R \times \mathcal{N}_0   \\
((u,h),(t,p)) &\mapsto (t+u,\rho(-t)(h) \cdot p)
\end{align}
where $\cdot$ is the action of $G_0$ on $\{0\}\times \mathcal{N}_0$. Indeed, denote this action by $*$ and write for $h \in G_0$, $h*(t,p) = h  \circ \phi^t(0,p)=\phi^t \circ \phi^{-t} \circ h \circ \phi^t (0,p)=\phi^t \circ \rho(-t)(h) (0,p)=\phi^t(0,\rho(-t)(h) \cdot p)=(t,\rho(-t)(h) \cdot p)$. 

\subsection{\textbf{$\dim \overline{\pi}(\Bar{\mathfrak{g}})=2:$ case of plane waves}}

Let $(M,g,V)$ be a locally homogeneous plane wave, and $\Bar{\mathfrak{g}}$ the Lie algebra of Killing fields preserving the line field $\R V$. Then, as in Paragraph \ref{Subsection 3.2}, we have a representation $\overline{\pi}: \Bar{\mathfrak{g}} \to \aff(\R)$ (see (\ref{Eq: representation of Lie algebra in aff(R)})).
In \cite[Proposition 4.3]{Leis}, the authors show that if $\dim \overline{\pi}(\Bar{\mathfrak{g}}) =2$, then $M$ is decomposable. In fact, it appears from their proof that under this assumption, $M$ is actually flat. To do so, they use the Killing field equations in local Brinkmann coordinates. 
In this paragraph, we give a coordinate-free proof of this fact. The key idea is that the existence of  a ``boost" in the isotropy of a locally homogeneous pp-wave implies that the metric is flat.
\begin{lemma}\label{Lemma: existence of a boost implies flatness}
Let $(M,g,V)$ be a locally homogeneous pp-wave. Let $p \in M$, and $I=\Stab(p)$. Let $\alpha: I \to \GL(T_p M)$ be the isotropy representation. If $\alpha(I)$ contains an element of the form 
	\begin{equation}
	\begin{pmatrix}
		\lambda & 0 & 0  \\
		0 & \lambda^{-1}  & 0  \\
		0 & 0 & A
	\end{pmatrix}
\end{equation}
with $\lambda^2 \neq 1$ and $A \in \O(n)$, then $M$ is flat. 
\end{lemma}
\begin{proof}
This forces all sectional curvatures to vanish at $p$, and hence everywhere  by local homogeneity. 

Since the leaves of $\mathcal{F}$ are flat and $V$ is parallel, it is sufficient to show that $R(U,X_i,U,X_i)=0$ for all $i$.  We have 
$$R(U,X_i,U,X_l) = \lambda^2 \sum_{j,k} A_{ji} A_{kl} R(U,X_j,U,X_k).$$
Define a matrix $S$ whose entries are $S_{jk}=R(U,X_j,U,X_k)$; it is a symmetric matrix. Then the above equality reads $S=\lambda^2 A{}^\top S A$. We also have $S^k=\lambda^{2k} A{}^\top S^k A$ for all $k \in \N$. Since $\lambda^2 \neq 1$, we have $\det(S) =0$, so that the kernel of $S$ is non-trivial. We claim that $S$ is nilpotent. Indeed, either $S^m=0$ for some $m \in \N$, or $\ker S^k = \ker S^{k+1}$ for some $k \in \N$. In the latter case, $\mathrm{Im}(S^k)$ is supplementary to $\ker S^k$, and both subspaces are preserved by $A$.  This allows to reduce the problem to a smaller dimension, and gives by induction that $S$ is nilpotent. Since $S$ is also symmetric, we must have $S=0$. 
\end{proof}
\begin{corollary}\label{Cor: dim rho(G)= 2 implies flat plane wave}
	Let $(M,g,V)$ be a locally homogeneous plane wave, and $\Bar{\mathfrak{g}}$ the Lie algebra of its Killing fields preserving the line field $\R V$. Consider the representation $\overline{\pi}: \Bar{\mathfrak{g}} \to \aff(\R)$. If $\dim \overline{\pi}(\Bar{\mathfrak{g}}) =2$, then $M$ is flat.	
\end{corollary}
	\begin{proof}
	If $\dim \overline{\pi}(\Bar{\mathfrak{g}}) =2$, there exists $K \in \Bar{\mathfrak{g}}$  whose local flow fixes the leaf $F(p)$ through $p$, and acts by homotheties on the space of leaves $M/\mathcal{F}$. Let $f$ be a local isometry of the flow. The induced action $d_p f : T_p(M/\mathcal{F}) \to T_p(M/\mathcal{F})$ has eigenvalue $\lambda \neq  1$, for otherwise $f$ would have trivial action on the space of leaves. Remember that $M$ admits an infinitesimal action of the Heisenberg algebra, that acts trivially on $M/\mathcal{F}$. We can compose $f$ by an element in $\heis_{2n+1}$ to get a new local isometry $f$ that fixes $p$, without altering the eigenvalue $\lambda$. So now we also suppose $f(p)=p$. This $\lambda$ is also an eigenvalue of $f$, i.e. there exists $U \in T_p M$ transverse to $V^{\perp}$ such that  $d_p f(U) = \lambda U$. We can suppose $g(U,V)=1$. Then necessarily, $V$ is an eigenvector with eigenvalue $\lambda^{-1}$. Define $E= \Vect(U,V)^{\perp}$: it is a Riemannian subspace of $T_p M$ preserved by $d_p f$. Choose an orthonormal basis $X_1,\dots,X_n$ of $E$: the matrix of $d_p f$ in the basis $\{U,V,X_1,\dots,X_n\}$ has the form:
		\begin{equation}
			\begin{pmatrix}
				\lambda & 0 & 0  \\
				0 & \lambda^{-1}  & 0  \\
				0 & 0 & A
			\end{pmatrix}
		\end{equation}
		with $A \in \O(n)$. We conclude using Lemma \ref{Lemma: existence of a boost implies flatness}.
	\end{proof}
\begin{remark}
	This proof doesn't work in the case of a locally homogeneous pp-wave which is not a plane wave. Indeed, since $\heis_{2n+1}$ does not act locally transitively on a leaf of $\F$, one does not obtain a boost in the isotropy as was done in the previous proof.  In \cite[Example 4.1]{Leis}, the authors provide an example of a $3$-dimensional homogeneous pp-wave, where $\dim \overline{\pi}(\Bar{\mathfrak{g}}) =2$,  which is not flat. 
\end{remark}
\section{Isometries with trivial action on the space of leaves}\label{Section: completeness of F-leaves}

Let $(\mathbf{X},g,V)$ be a $1$-connected homogeneous plane wave of dimension $n+2$. In this section, we begin the investigation of the Lie subgroup of $\Isom(\mathbf{X}, \R V)$ acting trivially on the space of leaves. This amounts to studying the completeness of Killing fields tangent to the $\F$-foliation. A first step is then to study the (geodesic) completeness of the $\F$-leaves of such manifolds. 
In particular, we show that the infinitesimal action of the Heisenberg algebra in the non-flat case integrates into a global isometric action of the Heisenberg group $\Heis_{2n+1}$. 

\subsection{Completeness along the $\F$-leaves}
Even if we are interested in the isometry group of $1$-connected homogeneous plane waves, we also deal with homogeneous pp-waves in the study of completeness of $\F$-leaves in this paragraph.

Let $(\mathbf{X},g,V)$ be a pp-wave of dimension $n+2$. Recall that the $\F$-leaves have a unimodular affine lightlike geometry (see Definition \ref{Introduction-Definition: Affine unimodular lightlike group} and Paragraph \ref{Section: Plane waves as curves in Aut(Heis)} ``Plane wave case''), i.e. a $(\mathsf{L}_\u(n), \R^{n+1})$-structure (in the sense of geometric structures, see \cite[Chapter 3]{thurston2022geometry}), where 
$\mathsf{L}_\u(n)=\O(n) \ltimes \Heis_{2n+1}$. 
\begin{definition}
Let $\mathbf{X}$ be a $1$-connected manifold with a $(G,X)$-structure, and $d: \mathbf{X} \to X$ a developing map. We say that $\mathbf{X}$ is homogeneous under the action of a group $L$ preserving the structure if there is a morphism $\theta: L \to G$ such that $d$ is $\theta$-equivariant.   
\end{definition}
If $\mathbf{X}$ is a $1$-connected non-flat homogeneous plane wave, the $\F$-leaves are homogeneous under the action of a group preserving the $(\mathsf{L}_\u(n), \R^{n+1})$-structure (see the proof of Corollary \ref{Corollary: homog plane wave has complete F-leaves} for the details). On the other hand, if $\mathbf{X}$ is a pp-wave (possibly flat), the $\F$-leaves are not necessarily homogeneous under the action of a group preserving the $(\mathsf{L}_\u(n), \R^{n+1})$-structure, but in this case, they are homogeneous under the action of a group preserving the $(\mathsf{L}(n), \R^{n+1})$-structure.

\begin{proposition}\label{Proposition: L_u-homogeneous implies F-complete}
Let $F$ be an $n+1$ dimensional $1$-connected manifold with a $(\mathsf{L}_\u(n),\R^{n+1})$-structure. If $F$ is homogeneous under the action of a group preserving the structure, then it is complete. 
\end{proposition}
\begin{proposition}\label{Proposition: L-homogeneous implies partial F-complete}
Let $F$ be an $n+1$ dimensional $1$-connected manifold, having a $(\mathsf{L}(n), \R^{n+1})$-structure. Then $F$ is naturally equipped with a Riemannian degenerate metric $g$. If $F$ is homogeneous under the action of a group preserving the structure, then it is either complete, or the incompleteness occurs along the $1$-dimensional  foliation defined by the kernel of $g$. 
\end{proposition}
The following lemma is needed in the proof of Propositions \ref{Proposition:   L_u-homogeneous implies F-complete} and \ref{Proposition: L-homogeneous implies partial F-complete}.

\begin{lemma}\label{Lemma: Homogeneous implies covering}
Let $M$ and $N$ be two smooth manifolds, with a local diffoemorphism $d: M \to N$. Assume that $M$ is $G_1$-homogeneous and $N$ is $G_2$-homogeneous, and there is a morphism $\theta: G_1 \to G_2$ such that $d$ is $\theta$-equivariant. Then the map $d$ is a covering from $M$ to $d(M)$. 
\end{lemma}
\begin{proof}
Let $\gamma: \R \to d(M)$ be a (continuous) curve, and let $x_0 \in M$ such that $d(x_0)=\gamma(0)$. We will show that there is a unique curve (lift) $\Tilde{\gamma}:\R \to M$ such that $d \circ \Tilde{\gamma}=\gamma$. We have that $d(M)$ is $\theta(G_1)$-homogeneous, so $\gamma(t)=g_t\gamma(0)$, where $g_t$ is a curve in $\theta(G_1)$. We claim that $g_t$ can be chosen continuously. To see this, set for $t \in \R$, $F_t:=\{g\in \theta(G_1), \ g\gamma(0)=\gamma(t)\}$.  
Set $G:=\theta(G_1)$, and consider the subgroup $H:=\Stab(\gamma(0))$. Define the map $\psi: \R \to G/H, \;\psi(t)=[g_t]  $, where $g_t\in F_t$. This map is continuous. Consider the pullback of the bundle map $\pi:G \to G/H$ by $\psi$. It has fibers $F_t$.
Since $\R$ is contractible this bundle has a global section, and the claim follows. Now, $\theta:G_1\to \theta(G_1)$ is a covering map ($G_1/\mathsf{\ker \theta}\cong \theta(G_1)$, where the kernel is discrete), so there is a unique lift (at the identity) of $g_t$ to a continuous curve $l_t$ in $G_1$, such that $\theta(l_t)=g_t$. Using the equivariance, we can write $\gamma(t)=g_t\gamma(0)=\theta(l_t)d(x_0)=d(l_t(x_0))$. Hence the lift of $\gamma$ is given by the curve $\Tilde{\gamma}(t)=l_t(x_0)$. 
\end{proof}

\begin{proof}[Proof of Proposition \ref{Proposition: L_u-homogeneous implies F-complete}] 
We have a developing map $d : F \to \R^{n+1}$ which is $\theta$-equivariant, where $\theta : G_0 \to \mathsf{L}_\u(n)$ is a  representation, and $G_0$ acts transitively on $F$. By Lemma \ref{Lemma: Homogeneous implies covering}, $d$ is covering map from $F$ to $d(F)\subset \R^{n+1}$, so to prove Proposition \ref{Proposition: L_u-homogeneous implies F-complete}, it is enough to show that $d$ is surjective onto $\R^{n+1}$. We have $\mathsf{L}_\u(n)= \Aut(\R^{n+1}, g_0=\sum_{i=1}^{n} dx_i^2, \partial_{x_0})$. Let $g$ (resp. $V$) be the pullback of $g_0$ (resp. $\partial_{x_0}$) by $d$. Let $ \mathcal{V}$ be the foliation defined by $V$.

The space of leaves $F/\mathcal{V}$ has a (Hausdorff) manifold structure. Indeed, it is homogeneous under the action of $G_0$,  and the stabilizer of any $\mathcal{V}$-leaf is a closed subgroup of $G_0$. The latter follows from the fact that the leaves are closed: a $V$-leaf is a connected component of the inverse image of a $\partial_{x_0}$ (closed) leaf in $\R^{n+1}$ by $d$, hence closed.  So $F/\mathcal{V}$ identifies with $G_0/I$, where $I$ is a closed subgroup of $G_0$. We claim that the space ${F/\mathcal{V}}$ is mapped bijectively to the space of $\partial_{x_0}$-leaves of $\R^{n+1}$. Indeed, the induced local diffeomorphism $\Bar{d}:=d_{|_{F/\mathcal{V}}}$ is a local isometry from $F/\mathcal{V}$, equipped with the Riemannian metric induced by $g$, to $(\R^n,g_0)$. And since ${F/\mathcal{V}}$ is homogeneous, it is complete, so the map $\Bar{d}$ is a covering, hence a diffeomorphism. On the other hand, the restriction of $d$ to a leaf of $V$ is a global diffeomorphism into a leaf of $\partial_{x_0}$. We claim that $d$ maps this leaf surjectively to a $\partial_{x_0}$-leaf. To see this, let $D_0$ be the leaf of $\partial_{x_0}$ through the origin, and let $H= \Stab(D_0)$ in $\mathsf{L}_\u(n)$. We have $H=(\O(n) \ltimes \R^n) \ltimes \R$, where $\O(n) \ltimes \R^n$ acts trivially on $D_0$ and $\R$ acts by translation. Let $\mathcal{V}_0$ be a leaf of $V$ whose image by $d$ is in $D_0$ (note that $D_0$ is indeed in the image of $d$ by the fact that $\Bar{d}$ is a diffeomorphism onto $\R^n$). Since $F$ is homogeneous, $\Stab(\mathcal{V}_0)$ acts transitively on $\mathcal{V}_0$, so  $\theta(\Stab(\mathcal{V}_0)) \subset H$ acts transitively on $\theta(\mathcal{V}_0) \subset D_0$. But a subgroup of $H$  either acts trivially on $D_0$ or has transitive action. We deduce that $\theta(\mathcal{V}_0)$ must be $D_0$. Now, if $D$ is any leaf of $\partial_{x_0}$, the $\Stab(D)$-action on $D$ is conjugate to the $H$-action on $D_0$,  hence the same conclusion holds for the $V$-leaves mapped to $D$.
\end{proof}
\begin{proof}[Proof of Proposition \ref{Proposition: L-homogeneous implies partial F-complete}] 
We have a developing map $d: F \to \R^{n+1}$ which is $\theta$-equivariant, where $\theta: G_0 \to \mathsf{L}(n)$ is a representation, and $G_0$ acts transitively on $F$. We have $\mathsf{L}(n)= \Aut(\R^{n+1}, g_0=\sum_{i=1}^{n} dx_i^2, \R \partial_{x_0})$. Let $g$ (resp. the line field $l$) be the pullback of $g_0$ (resp. $\R \partial_{x_0}$) by $d$. We denote  by $\mathcal{L}$ the foliation defined by $l$. Here again, $d$ is a covering map from $F$ to $d(F)$ which induces a diffeomorphism from $F / \mathcal{L}$ to $\R^n := \R^{n+1}/ \R \partial_{x_0}$. Let now $D_0$ be the leaf of $\partial_{x_0}$ through the origin, and let $H:= \Stab(D_0)$ in $\mathsf{L}(n)$. We have $H= ((\R \times \O(n)) \ltimes \R^n) \ltimes \R$, where $\O(n) \ltimes \R^n$ acts trivially on $D_0$ and $\mathsf{Aff}(\R)$ acts by its usual action. Let $\mathcal{L}_0$ be a leaf of $\mathcal{L}$ whose image by $d$ is in $D_0$. Then $\theta(\Stab(\mathcal{L}_0)) \subset H$ acts transitively on $d(\mathcal{L}_0) \subset D_0$. And a subgroup of $H$ which does not act trivially either acts transitively, in which case $d(\mathcal{L}_0)=D_0$, or $d(\mathcal{L}_0)$ is half a line. The same conclusion holds for any leaf of $\mathcal{L}$, and they are moreover either all complete or all incomplete. Hence the proposition. 
\end{proof}
\begin{corollary}\label{Corollary: homog plane wave has complete F-leaves}
A homogeneous non-flat plane wave $(\mathbf{X},g,V)$ has complete $\mathcal{F}$-leaves. 
\end{corollary}
\begin{proof}
One may suppose that $\mathbf{X}$ is connected, and $1$-connected up to taking the universal cover. 
Let $G=\Isom^{\mathsf{o}}(\mathbf{X},\R V)$. Since $\mathbf{X}$ is non-flat, $V$ is the unique parallel lightlike vector field of $\mathbf{X}$ (up to scale), hence $\Isom(\mathbf{X},g)=\Isom(\mathbf{X},g, \R V)$, and then $\mathbf{X}$ is homogeneous under the action of $G$. 
The stabilizer of a leaf preserves the induced flat connection, the induced degenerate Riemannian scalar product together with  the lightlike direction $\R V$, hence represents in the affine lightlike group $\mathsf{L}(n)$.  On the other hand, we have $\dim \pi(G) = 1$ by Corollary \ref{Cor: dim rho(G)= 2 implies flat plane wave}. Due to homogeneity of $\mathbf{X}$, the $1$-parameter group with non trivial $\pi$-projection acts transversely on the $\F$-leaves. Therefore,  the subgroup of $G$ preserving a leaf of $\F$ is contained in the kernel of $\pi$, i.e. represents in the subgroup $\mathsf{L}_\u(n)$ of $\mathsf{L}(n)$. It follows that a leaf of $\F$ is homogeneous under the action of a group preserving the $(\mathsf{L}_\u(n),\R^{n+1})$-structure, hence the completeness by Proposition \ref{Proposition: L_u-homogeneous implies F-complete}.  
\end{proof}
\begin{corollary}\label{Corollary: homog pp-wave with V complete has complete F-leaves}
A homogeneous pp-wave (flat or non-flat) such that $V$ is complete, has complete $\mathcal{F}$-leaves. 
\end{corollary}
\begin{proof}
We Keep the same notation as in the previous proof. Here again, the stabilizer of a leaf of $\F$ in $G$ has a representation in $\mathsf{L}(n)$, and then a leaf of $\F$ is homogeneous under the action of a group preserving the $(\mathsf{L}(n), \R^{n+1})$-structure. The conclusion follows from Proposition \ref{Proposition: L-homogeneous implies partial F-complete}. 
\end{proof}
\begin{example}[Homogeneous pp-wave with incomplete $\F$-leaves]\cite[Example 7.113]{Besse}
Consider the Minkowski plane $\Mink^{1,1}:=(\R^{2},2 dudv)$. 
There is a subgroup $\mathsf{Aff}(\R)$ of  $\O(1,1) \ltimes \R^2$ acting transitively on the half Minkowski space $U=\{(v,u) \in \R^2, v >0\}$ as follows $(v,u) \mapsto (e^a v, e^{-a} u + b)$, where $(a, b) \in \R \ltimes \R$. Here, $V=\partial_v$ is a parallel lightlike vector field, so $U$ is a (homogeneous) pp-wave. And the $V^{\perp}$-leaves given by the $u$-levels are incomplete. 
\end{example}

\subsection{Global $\Heis_{2n+1}$-action}
As a consequence of Corollary \ref{Corollary: homog plane wave has complete F-leaves}, we obtain that for a $1$-connected non-flat homogeneous plane wave, the (global) Killing fields tangent to the $\F$-foliation are complete, hence the infinitesimal isometric action of the Lie algebra $\mathfrak{g}_0$ (see Equation (\ref{g_0})) integrates into a Lie group action of the universal cover of $G_0$. 

\begin{corollary}\label{Corollary: global Heis-action}
Let $\mathbf{X}$ be a $1$-connected non-flat homogeneous plane wave of dimension $n+2$. The (global) Killing fields tangent to the $\F$-foliation are complete. Hence $\mathbf{X}$  admits a $\Heis_{2n+1}$ global isometric action, preserving individually the $\F$-leaves, and acting transitively on each leaf of $\F$.     
\end{corollary}

\section{Isometry group of $1$-connected homogeneous plane waves}\label{Section: Isometry group}
Let $\mathbf{X}$ be a $1$-connected non-flat homogeneous plane wave. In this section we compute the identity component of the isometry group of $\mathbf{X}$ (Theorem \ref{THEOREM: Isometry group homog case}), and obtain a description of such spaces as $1$-parameter groups in $\Aut(\Heis_{2n+1})$. Note that in the non-flat case, $\mathbf{X}$ admits a unique parallel lightlike vector field $V$ (up to scale). Indeed, if we assume the existence of another parallel lightlike vector field $V_1$, then, the plane wave locally decomposes as an orthogonal product $P \times X_0$, where $P$ is the flat Lorentzian surface tangent to  span$(V, V_1)$, and $X_0$  is a Riemannian space orthogonal to $P$, and thus to $V$.   In particular, due to the curvature conditions on plane waves (and more generally, on pp-waves), $X_0$ must necessarily be flat. This implies that $\mathbf{X}$ is  also flat. It is worth mentioning that these curvature conditions exclude, for instance, the possibility of having the product of a Cahen-Wallach space (or a Minkowski space) with a non-flat Riemannian space as a pp-wave; such products instead fall within the broader class of Brinkmann manifolds.

So, for a $1$-connected non-flat homogeneous plane wave $\mathbf{X}$, we have $\Isom(\mathbf{X})=\Isom(\mathbf{X}, \R V)$. Define $\widehat{G}:= \mathsf{Isom}(\mathbf{X})$.\\

The notations used in this section for the subgroups of $\Heis_{2n+1}$ are given in the introduction (Subsection \ref{Subsection-Introduction: deformation of Heis action}). Let us recall them here. The standard action of $\Heis_{2n+1}=\R^n \ltimes \R^{n+1}$ on $\R^{n+1}$ is the one in which $A^-\times Z=\R^{n+1}$ acts by translation and $A^+= \R^{n}$ acts linearly by  unipotent matrices.  We use the usual decomposition $\heis_{2n+1} = \a^+ \oplus \a^- \oplus \z$ as a vector space, where $\z$ is the center of $\heis_{2n+1}$, $\a^+$ and $\a^-$ the Lie algebras of $A^+$ and $A^-$ respectively. 
Denote by $X_1,\dots,X_n$ (resp. $Y_1,\dots,Y_n$) a basis of $\a^+$ (resp. $\a^-$), such that all Lie brackets are zero but $[X_i,Y_i]=z$, for $i=1,\dots,n$.

\subsection{Isometry group}
Throughout this section, $\mathbf{X}$ is a $1$-connected non-flat homogeneous plane wave of dimension $n+2$. In the following we give the form of the identity component of the isometry group of $\mathbf{X}$.

\begin{proposition}\label{Proposition: isometry group1} Let $\mathbf{X}$ be a $1$-connected non-flat homogeneous plane wave. Let $\widehat{G}= \mathsf{Isom}(\mathbf{X})$ and $G$ its identity component. Then 
\begin{itemize}
    \item[(1)] The index $[\widehat{G}:G]$ is finite.
    \item[(2)] $G \cong (\R\times K)\ltimes \Heis_{2n+1}$,  where $K$ is some closed connected subgroup of $\O(n)$ acting on $\Heis_{2n+1}$, trivially on the center and by the standard action on $A^+$ and $A^-$.
\end{itemize}
And $\mathbf{X}$ identifies with the homogeneous space $(\R\times K)\ltimes \Heis_{2n+1}/ K \ltimes A^+$. \medskip
\end{proposition}
The $\R$-action on $\Heis_{2n+1}$ and the compact group $K$ will be characterized in Theorem \ref{THEOREM: Isometry group homog case}.  \medskip

Before proving Proposition \ref{Proposition: isometry group1}, let us make some observations. Remember that we have an exact sequence (see ($\ref{Exact sequence on G}$) in Section $\ref{Section: Affine structure on the space of leaves}$):
$$1 \rightarrow  G_0 \rightarrow G \rightarrow  \pi(G) < \mathsf{Aff}(\R) \rightarrow 1. $$

Since $\mathbf{X}$ is non-flat, the projection to the affine group is one dimensional (Corollary \ref{Cor:  dim rho(G)= 2 implies flat plane wave}). Hence, by Proposition \ref{Propositiion: Isometry group of Homogeneous Brinkman}, we get $G$ as a semi-direct product  $\R\ltimes G_0$, where $G_0$ is a Lie group with a  (necessarily faithful) representation $\theta: G_0 \to \mathsf{L}_\u(n)= \O(n) \ltimes \Heis_{2n+1}$,  and $\R$ acts on it through a $1$-parameter group in $\Aut(G_0)$. So $G$ is isomorphic to  $\R \ltimes \theta(G_0)$, where the $\R$ action on $\theta(G_0)$ is the conjugation of the $\R$ action on $G_0$ by $\theta:G_0 \to \theta(G_0)$. Moreover, by Corollary \ref{Corollary: global Heis-action}, $G_0$ contains a subgroup $H$ isomorphic to the Heisenberg group. By Lemma \ref{Lemma: rho(Heis)=Heis}, we have that $\theta(H)= \Heis_{2n+1}$.
It follows in particular that $\theta(G_0)$ contains $\Heis_{2n+1}$, 
hence it itself splits as $K\ltimes \Heis_{2n+1}$, where $K := \theta(G_0) \cap \O(n)$ is some subgroup of $\O(n)$. Finally, $G$ is isomorphic to $\R \ltimes (K \ltimes \Heis_{2n+1})$. 

Next, we identify $\mathbf{X}$ as a homogeneous space of the Lie group $\R \ltimes (K \ltimes \Heis_{2n+1})$.  
In Section \ref{Section: Affine structure on the space of leaves}, we wrote a diffeomorphism $\mathbf{X} \cong \R \times \mathcal{N}_0$, where $\mathcal{N}_0$ is a fixed leaf of $\F$ (which is necessarily $1$-connected). The $G$ action on $\mathbf{X}$ is then equivalent to a $\R \ltimes G_0$ action on $\R \times \mathcal{N}_0$ given by (\ref{Eq: R ltimes G_0 action on RxN_0}). 
Now, by the previous section, we know that the $\F$-leaves are complete, so we have a developing map $d: \mathcal{N}_0 \to \R^{n+1}$ which is $\theta$-equivariant, where $\theta$ is the representation given above. Thus, for $g \in G_0$, we have $\theta(g) * d(p) = d(g \cdot p)$, where $*$ is the standard action of $\theta(G_0)=K\ltimes \Heis_{2n+1} \subset \mathsf{L}_\u(n)=\O(n)\ltimes \Heis_{2n+1}$ on $\R^{n+1}$. And the stabilizer at the origin for this action is $K \ltimes A^+$. So $\mathbf{X}$ identifies with $\R \ltimes (K \ltimes \Heis_{2n+1})/ K \ltimes A^+$.

\begin{observation}\label{Observation: L preserves heis decomposition}
Since $\Heis_{2n+1}$ is the nilradical of $\theta(G_0)$, it is preserved by any automorphism of $\theta(G_0)$. Hence the $\R$-factor normalizes $\Heis_{2n+1} \subset \theta(G_0)$, and $G$ contains the subgroup $\R \ltimes_{\rho} \Heis_{2n+1}$, where $\rho: \R \to \Aut(\Heis_{2n+1}), \rho(t) = e^{t L}$, and $L \in \Der(\Heis_{2n+1})$. Replacing $L$ by $L+\ad_h$ for a suitable $h \in \heis_{2n+1}$ (which amounts to taking another $1$-parameter group in $G$ as a splitting of the exact sequence above), we can suppose that $L$ preserves the decomposition $\heis_{2n+1}=\C^n \oplus \z$ as a vector space. 
\end{observation}

To get the form of $G$ given in Proposition \ref{Proposition: isometry group1}, we need the following lemma.
    \begin{lemma}\label{product structure}
    Let $Q$ be the semi-direct product $Q=\R\ltimes K$ where $K$ is a compact connected Lie group. Then $Q$ is isomorphic to the product $\R\times K$. 
\end{lemma}
\begin{proof} 
Let $\c$ be the Lie algebra of $K$. Then $\c$ has a unique decomposition as the sum of two ideals, one of which is a (maximal) abelian ideal $\r$ and the other is a semi-simple ideal $\s$. To see this, consider the Levi-decomposition $\c=\s\oplus\r$, where $\r$ is the radical of $\c$ and $\s$ is a Levi factor. The Lie subgroup $R$ with Lie algebra $\r$ is closed in $K$, for otherwise one can consider its closure, whose Lie algebra is then a  solvable Lie algebra strictly bigger than $\r$, contradicting the fact that $\r$ is maximal.  As a consequence, $\r$ is a solvable Lie algebra of a compact group, hence necessarily abelian, and $R$ is isomorphic to a torus $\T^n$. The uniqueness amounts to proving that $\s$ is unique. Let $\s'$ be  another Levi factor. Consider the projection $\c=\s \oplus \r \to \r$. Since $\s'$ is semi-simple and $\r$ is abelian, the image of $\s'$ by the projection is trivial, hence contained in the kernel $\s$. This implies $\s'=\s$. Now, let $\widetilde{K}$ be the universal cover of $K$. It follows from the previous discussion that $\widetilde{K}$ is isomorphic to $\R^n\ltimes \widetilde{S}$, where $\widetilde{S}$ is semi-simple. The $\R$ action on $K$ is given by some $1$-parameter group in $\Aut(K)$, denote it by $\phi^t$, and let $\widetilde{\phi}^t$ be its lift to $\widetilde{K}$. 
First, it follows from the first part of the proof that for any $t\in \R$,
 $\widetilde{\phi}^t(\R^n)=\R^n$ and
 $\widetilde{\phi}^t(\widetilde{S})=\widetilde{S}$.
In particular, $\phi^t$ preserves $R$, which is a torus in $K$, and the induced action on it is trivial since the automorphism group of a torus is discrete. So $\widetilde{\phi}^t_{\vert \R^n}$ is also trivial. Now, since $\widetilde{S}$ is semi-simple, any derivation is inner, so we can assume that the $\widetilde{\phi}^t$-action on it is trivial up to conjugacy by a $1$-parameter subgroup of $\widetilde{S}$, say $f(t)$.  And we have  $f(t)\widetilde{\phi}^t f(t)^{-1}=\mathrm{Id_{\widetilde{K}}}$. Conjugation by $f(t)$ passes to the quotient (since $\pi_1(K)$  is in the center of $\widetilde{K}$), and allows to get a $1$-parameter group with trivial action on $K$.
\end{proof}

\begin{proof}[Proof of Proposition $\ref{Proposition: isometry group1}$]
$\mathsf{(1)}$   We can make a general argument as follows. Let $I$ be the full isotropy of a point in $\mathbf{X}$, and $I^\circ=I \cap G$ its identity component. Since $\mathbf{X}$ is homogeneous, and its Killing fields are  integrable, a standard fact from Gromov's rigid transformation groups theory \cite[Theorem 3.5.A]{Gro} states that the image of $I$ inside $\O(n+1,1)$ by the faithful representation is an algebraic group, hence has finitely many connected components.  It follows that $G$ has finite index in $\widehat{G}$ (in fact has the same index as $I^\circ$ in $I$).

$\mathsf{(2)}$ From the discussion after the proposition, we have $G\cong\R \ltimes (K \ltimes \Heis_{2n+1})$. We will first show that $K$ is  closed in $\O(n)$. Let $I$ be the isotropy of some point in $\mathbf{X}$, we have that $\theta(I)$ is algebraic (from Gromov's rigid transformation groups theory), hence closed in $\O(n) \ltimes \Heis_{2n+1}$ and decomposes as $\theta(I)= K' \ltimes A$, with $K'$ a subgroup of $\O(n)$ and $A \subset A^+$. In particular, $K'=(K' \ltimes A) \cap \O(n)$ is closed in $\O(n)$. But $K'=K$, and the claim follows.  

Now, we will show that the $\R$ action on $K \ltimes \Heis_{2n+1}$ preserves $K$. For this, it is sufficient to show that the $\ad$-action of $L$ preserves $\mathfrak{c}$. Let $A \in \mathfrak{c}$. It follows from the identity $\ad_{[A,L]}=[\ad_A, \ad_L]$ and the fact that $L$ and $\mathfrak{k}$ preserve the decomposition $\heis_{2n+1} = \C^n \oplus \R$ as a vector space that the $\ad$-action of  $[A,L]$ preserves the decomposition too. On the other hand, we have  $[A,L] \in \mathfrak{c} \oplus \heis_{2n+1}$, since $\mathfrak{c} \oplus \heis_{2n+1}$ is an ideal. By looking at the action of $\ad_U$ on $\heis_{2n+1}$, for $U \in \mathfrak{c} \oplus \heis_{2n+1}$, we see that necessarily $[A,L] \in \mathfrak{c} \oplus \z$. So the action of the $\R$-factor sends $K$ into $K \times Z$. But since $K$ is compact, the image must be in $K$, since an element $(c,z)$, with $z \neq 0$, generates an unbounded sequence in $K \times Z$ (here $K$ acts trivially on $Z$). So far we proved that $K$ is preserved by the $\R$-action. Now applying Lemma $\ref{product structure}$, we can replace the $\R$-action by changing the $1$-parameter group in the product $\R\ltimes K$, to have a product structure.
\end{proof}

\subsection{Homogeneous plane waves vs $1$-parameter groups in $\Aut(\Heis_{2n+1})$}\label{Section: Homogeneous plane waves vs 1-parameter groups in AUt(Heis)}
Let $(\mathbf{X},g,V)$ be a $1$-connected homogeneous non-flat plane wave of dimension $n+2$. 
We know from Proposition \ref{Proposition: isometry group1} that the connected component of the isometry group is isomorphic to $G_\rho=(\R\times K)\ltimes_{\rho} \Heis_{2n+1}$, where $K$ is some connected closed subgroup of $\SO(n)$, and that $\mathbf{X}$ identifies as $X_{\rho}=G_\rho/I$, with $I=K\ltimes A^+$.  This motivates the following question, whose answer is given in Theorem \ref{THEOREM: Isometry group homog case}.

\begin{question}\label{Question: existence of G_rho invariant metric}
Let $G_{\rho}=(\R\times K)\ltimes_{\rho} \Heis_{2n+1}$ such that 
\begin{itemize}
    \item $K$ is a connected compact group acting on $\Heis_{2n+1}$ trivially on the center, and by preserving the decomposition $A^+ \oplus A^-$,
    \item $\R$ acts on $\Heis_{2n+1}$ through a representation $\rho: \R \to \Aut(\Heis_{2n+1})$, $\rho(t)=e^{tL}$, with $L \in \Der(\heis_{2n+1})$. Suppose that $L$ preserves the decomposition $\heis_{2n+1}=\C \oplus \z$ as a vector space.  
\end{itemize}
Consider the homogeneous space $X_{\rho}=G_{\rho}/I$, where $I=K \ltimes A^+$. 
\begin{enumerate}
    \item Which spaces $X_{\rho}$ admit a Lorentzian metric $g$ invariant under the left action of $G_\rho$? We know that if such a metric exists, it is necessarily a plane wave metric (Proposition \ref{Proposition: action preserves Lorentz metric => plane wave}). 
    \item If such a metric exists, is it unique (up to isometry of $X_\rho$)?  
\end{enumerate}
\end{question}

\subsubsection{\textbf{$G_\rho$-invariant Lorentzian metrics on $X_\rho$ up to isometry}}
$G_\rho$ preserves a Lorentzian metric on $G_\rho/I$ if and only if the $\Ad(I)$-action on $\mathfrak{g}/\mathfrak{i}$ preserves a Lorentzian scalar product $q$. 
Equivalently, $\ad_h$ is skew-symmetric with respect to $q$ for any $h \in \mathfrak{i}$. Let $T$ be a basis of the $\R$-factor. We have $\ad_h(T)=L(h)$, for any $h \in \heis_{2n+1}$. We have  $\mathfrak{g}_\rho/\mathfrak{i} \cong \R \oplus \R \oplus \mathfrak{a}^-  $ with basis $(T, z, Y_1,\dots,Y_n)$. We write $q$ with respect to this basis. Denote by $S_\rho$ the set of all $\ad_{\mathfrak{i}}$-invariant Lorentzian inner products on $\mathfrak{g}_\rho/\mathfrak{i}.$ \\

\noindent \textbf{Action of $\Aut_{\mathfrak{i}}^1(\mathfrak{g}_\rho)$ on $S_\rho$.} Denote by $\Aut_I(G_\rho)$ (resp. $\Aut_{\mathfrak{i}}(\mathfrak{g}_\rho)$) the automorphism group of $G_\rho$ (resp. $\mathfrak{g}_\rho)$ preserving $I$ (resp. $\mathfrak{i}$).  Let $f \in \Aut_I(G_\rho)$, and denote by $\bar{f}$ the induced diffeomorphism on $X_\rho$. Let $q_0$ and $q_1 \in S_\rho$. And let $g_0$ and $g_1$ be the $G_\rho$ left invariant Lorentzian metrics on $X_\rho$ determined by $q_0$ and $q_1$ respectively. If $d_e f \in \Aut_{\mathfrak{i}}(\mathfrak{g})$ induces an isometry from $(\mathfrak{g}/\mathfrak{i}, q_0)$ to $(\mathfrak{g}/\mathfrak{i}, q_1)$, then $\Bar{f} : (X_\rho,g_0) \to (X_\rho, g_1)$ is an isometry. Let $\Aut_{\mathfrak{i}}^1(\mathfrak{g}_\rho) \subset \Aut_{\mathfrak{i}}(\mathfrak{g}_\rho)$ be the subgroup of elements that are a lift of a Lie group isomorphism. Then elements in $S_{\rho}$ in the same $\Aut_{\mathfrak{i}}^1(\mathfrak{g}_\rho)$-orbit (for the pullback action) give rise to isometric metrics on $X_\rho$.\\

Let $\mathfrak{g}_\rho= (\R \oplus \mathfrak{k}) \ltimes \heis_{2n+1}$. As mentioned in Observation \ref{Observation: L preserves heis decomposition}, we can suppose without loss of generality that $L$ preserves the decomposition $\heis_{2n+1}=\C^n \oplus \z$ as a vector space. 
Hence  
\begin{equation}\label{Eq: Derivation matrix L}
		L= \begin{pmatrix}
			A & B & 0  \\
			D & C & 0  \\
			0 & 0 & \delta 
		\end{pmatrix} 
\end{equation}
is the matrix representing $L$ in $\mathfrak{a}^+ \oplus \mathfrak{a}^- \oplus \mathfrak{z}$.

\begin{lemma}\label{Lemma: derivations of heis}
The matrix $L$ above is a derivation of $\heis_{2n+1}$ if and only if $B$ and $D$ are symmetric, and $C+A=\delta I_n$, where $L(z)= \delta z$. 
 \end{lemma}
 \begin{proof}
  This is a straightforward computation.    
 \end{proof}

The following proposition gives a necessary condition on the $\R$-action for the existence of $g$. 
\begin{proposition}\label{Proposition: G preserves metric on G/I implies D definite}
If $X_\rho=G_\rho /I$ admits a $G_\rho$-invariant Lorentzian metric, then $D$ is a definite symmetric matrix. In particular, this implies that $L(\mathfrak{a}^+) \cap (\mathfrak{a^+} \oplus \mathfrak{z}) = 0$. 
\end{proposition}

The following well known lemma will be used in the proof of Proposition \ref{Proposition: G preserves metric on G/I implies D definite} above. We didn't find a reference for the proof, so we give a proof here. 
\begin{lemma}\label{Lemma: f nil degree 3}
A nilpotent non-zero endomorphism $f$ of a vector space $E$ is an infinitesimal isometry of some Lorentzian scalar product $q$ if and only if its nilpotency order equals $3$, and $E=E_1 \oplus E_2$, where $E_1$ is a $3$-dimensional subspace on which $f$ is $3$-nilpotent, and $E_2$ is $q$-orthogonal to $E_1$, with $f_{\vert E_2}=0$. 
\end{lemma}
\begin{proof}
Let $f$ be nilpotent of degree $3$. By Jordan decomposition, $f$ is conjugate to a matrix containing $3$-nilpotent  blocs $J_i$ of dimension $3$, with $J_i(e_0^i)=0, J_i(e_1^i)=e_0^i$ and $J_i(e_2^i)=e_1^i$. Suppose that $f$ is skew-symmetric for some quadratic form $q$.  Then  $q(f(e_1^i), e_0^j) + q(e_1^i, f(e_0^j))=0$. This yields $q(e_0^i,e_0^i)=q(e_0^j,e_0^j)=0$, and $q(e_0^i,e_0^j)=0$. For $q$ to be of Lorentzian signature, we must then have $i=j$, hence only one $3$-nilpotent bloc $J$. It follows that $f$ acts on a $3$-dimensional subspace $E_1$ by a $3$-nilpotent matrix $J$ as above, and is zero on a supplementary $E_2$. Then define $q$ such that  $q(e_0,e_0)=0$, $q(e_0,e_1)=0$, $q(e_2,e_1)=0$, $E_2$ is $q$-orthogonal to $E_1$, and $q_{E_2}$ is any Euclidean quadratic form. Then $q$ is a Lorentzian quadratic form, and $f$ is skew-symmetric for $q$. Conversely, if $f$ is nilpotent of degree $n \neq 3$, it is easy to see that any quadratic form $q$ for which $f$ is skew-symmetric has  a totally isotropic subspace of dimension $\geq 2$.
\end{proof}
\begin{proof}[Proof of Proposition \ref{Proposition: G preserves metric on G/I implies D definite}]
Remember that $G$ preserves a Lorentzian metric on $G/I$ if and only if the $\Ad(I)$-action on $\mathfrak{g}/\mathfrak{i}$ preserves a Lorentzian scalar product $q$. In particular, $\ad_h$ is skew-symmetric with respect to $q$, for any $h \in \mathfrak{a}^+$. For any $X \in \heis_{2n+1}$, $\ad_h^2(X)=0$, and $\ad_h(T) \in \heis_{2n+1}$. So $\ad_h^3=0$. The action of $\ad_h$ on $\mathfrak{g}/\mathfrak{i}$ has exactly degree $3$ if and only if $\ad_h^2(T)=[L(h),h] \neq 0$, for any $h \in \mathfrak{a}^+$.  	But $$[L(h),h] = [D(h),h] = -\langle D(h), h \rangle z. $$ And this is equivalent to $D$ being definite. By Lemma \ref{Lemma: f nil degree 3}, this gives a necessary condition for the existence of $q$.    
\end{proof}

\textbf{Reductions:} Let $G_{\rho_1}=(\R \times K)\ltimes_{\rho_1} \Heis_{2n+1}$ and $G_{\rho_2}=(\R \times K)\ltimes_{\rho_2} \Heis_{2n+1}$ be two isomorphic split extensions of $\R \times K$ by $\Heis_{2n+1}$. Let $f: G_{\rho_1} \to G_{\rho_2}$ be the isomorphism between the two extensions, and suppose that $f(I)=I$. Let $q_1 \in S_{\rho_1}$ and $q_2 \in S_{\rho_2}$, and suppose that $f$ induces an isometry from $(\mathfrak{g}_1/\mathfrak{i}, q_1)$ to $(\mathfrak{g}_2/\mathfrak{i}, q_2)$. Then $f: (X_{\rho_1}, g_1) \to (X_{\rho_2}, g_2)$ is an isometry, where $g_1$ (resp. $g_2$) is the left invariant Lorentzian metric on $X_{\rho_1}$ (resp. $X_{\rho_2}$) determined by $q_1$  (resp. $q_2$). This allows to make reductions of type: ``up to considering an isomorphic split extension, we can suppose that the matrix derivation (\ref{Eq: Derivation matrix L}) has a given shape".
Let us introduce a special case of extension isomorphism that will be used repeatedly in the sequel. Let $\phi \in \Aut(\Heis_{2n+1})$ preserving $A^+$, and consider  
\begin{align*}
        (\R \times K)\ltimes_\rho \Heis_{2n+1} &\to (\R \times K) \ltimes_{\rho'} \Heis_{2n+1}\\ 
        (t,k,h) &\mapsto (t,k,\phi(h))
    \end{align*}
where in $G_{\rho'}$, both $\R$-action and $K$-action are conjugated by the automorphism $\phi$ of $\Heis_{2n+1}$.
This is a Lie group isomorphism  preserving $I$.\\

\noindent \textbf{Reduction 1.} As a consequence of Proposition \ref{Proposition: G preserves metric on G/I implies D definite} above,  when dealing with Question \ref{Question: existence of G_rho invariant metric}, one can reduce to $G_\rho$ for which the representation $\rho(t)=e^{t L}$ is given by a derivation of the form 
\begin{equation}\label{Eq: New derivation L}
 L= \begin{pmatrix}
    A & B & 0\\
    I_n & \delta I_n - A^\top & 0 \\
    0 & 0 & \delta
\end{pmatrix},   
\end{equation}
with $\delta \in \R$. Let us justify this fact. Let $G_\rho$ such that the matrix $D$ in the derivation $L$ is definite. Since $D$ is symmetric definite, there exists $P \in \O(n)$ such that $P DP^\top =\mathrm{Diag}(\lambda_1,\dots,\lambda_n)$, where the coefficients  $\lambda_i \neq 0$ are either all positive or all negative. Define also $\Lambda:=\mathrm{Diag}(\sqrt{\vert \lambda_1 \vert},\dots,\sqrt{\vert \lambda_n \vert})$, and finally $P_{\Lambda}=\Lambda P \in \GL_n(\R)$. 
The matrix 
$\begin{pmatrix}
    P_{\Lambda} & 0\\
    0 & (P_{\Lambda}^\top)^{-1}
\end{pmatrix} \in \mathsf{Sp}_{2n}(\R)$, so it lifts to a (Lie algebra) automorphism $J= \begin{pmatrix}
    P_{\Lambda} & 0 & 0\\
    0 & (P_{\Lambda}^\top)^{-1} & 0\\
    0 & 0 & 1
\end{pmatrix}$ of $\heis_{2n+1}$ preserving $\a^+$. Let also 
$J_0=\begin{pmatrix}
    I_n & 0 & 0\\
    0 & -I_n & 0\\
    0 & 0 & -1
\end{pmatrix} \in \Aut_{\a^+}(\heis_{2n+1})$. 
Then $J_0J L (J_0J)^{-1}=
\begin{pmatrix}
    A' & B' & 0\\
    I_n & C' & 0\\
    0 & 0 & \delta
\end{pmatrix}$. From this we deduce that up to isomorphism, when $D$ is definite, the derivation $L$ represented in $\a^+ \oplus \a^- \oplus \z$ in (\ref{Eq: Derivation matrix L}) can be taken with $D=I_n$. \\

Recall that $\mathfrak{g}_\rho/\mathfrak{i} \cong \R \oplus \R \oplus\mathfrak{a}^-  $ with basis $(T, z, Y_1,\dots,Y_n)$. We write $q$ with respect to this basis. 
\begin{proposition}\label{Proposition: when G preserves a metric on G/I}
Let $G_\rho=(\R \times K) \ltimes \Heis_{2n+1}$ with $\rho(t)=e^{t L}$ as in (\ref{Eq: New derivation L}), and let $I=K \ltimes A^+$. 
Then $X_\rho=G_\rho /I$ admits a $G_\rho$-invariant Lorentzian metric if and only if $K$ is the (connected) Lie subgroup of $\SO(n)$ with Lie algebra $\mathfrak{k}=\{F \in \mathfrak{o}(n), \;[F, A]=[F,B]=0 \}$. It acts on $\Heis_{2n+1}$ trivially on the center and by its standard action on $A^+$ and $A^-$. 
Moreover, this metric is unique (up to isometry of $X_\rho$).
\end{proposition}
\begin{proof}
\textbf{Existence:}   
Writing the fact that $\ad_h$ is skew-symmetric for a quadratic form $q$ on $\mathfrak{g}/\mathfrak{i}$ for any $h \in \a^+$ yields :
\begin{itemize}
      \item $q(z,z)=0$
      \item $q(z, Y_i)=0$ for all $i=1,\dots,n$
      \item $q(T, Y_i)=0$ for all $i=1,\dots,n$
      \item $H=- q(z,T) I_n$, where $H=(q(Y_i,Y_j))_{i,j}$
\end{itemize}
If $q$ has Lorentzian signature, then necessarily $q(z,T) <0$. Moreover, we see that $q$ is completely determined, up to choosing $q(z,T)$ and $q(T,T)$. On the other hand, $\ad_h$, for $h \in \mathfrak{k}$, acts on $\heis_{2n+1}$ via the derivation matrix 
$\begin{pmatrix}
 -F^\top & 0 & 0\\
 0 & F & 0 \\
 0 & 0 & 0
\end{pmatrix}$, with $F \in M_n(\R)$. 
A straightforward computation shows that the $\R$ action and $K$ action commute if and only if $[F,A]=[F,B]=0$. 
Now we write the fact that $\ad_h$ is skewsymmetric for $q$. Using the scalar products already determined above, and the fact that $\ad_h(T)=0$, the only non-trivial condition is $q(\ad_h(Y_i),Y_j)+q(Y_i,\ad_h(Y_j))=0$ for any $i,j$. This means that $F$ is skewsymmetric for the Euclidean scalar product, hence an element of $\mathfrak{o}(n)$. If follows that $X_\rho$ admits a $G_\rho$-invariant Lorentzian metric if and only if $\mathfrak{k}=\{F \in \mathfrak{o}(n), [F,A]=[F,B]=0 \}$. \\
\textbf{Uniqueness:} To finish the proof, we will show that all such scalar products $q$ on $\mathfrak{g}_\rho /\mathfrak{i}$ define isometric metrics on $G_\rho / I$. Fix such $q$. And consider the Lie algebra isomorphism 
\begin{align*}
    \phi:\; &\R  \ltimes \heis_{2n+1} \to \R  \ltimes \heis_{2n+1} \\
    &\phi(T)=T + \frac{1}{2} \alpha z, \\
    &\phi(z)=\vert \beta \vert z,\; \phi(Y_i)= \sqrt{\vert \beta \vert} Y_i,\; \phi(X_i) = \sqrt{\vert \beta \vert} X_i, 
\end{align*}
where $\beta:=q(T,z)$ and $\alpha:=q(T,T)$. It preserves $\mathfrak{a}^+$, and induces an isometry from $(\R  \ltimes \heis_{2n+1}/ \mathfrak{a}^+, q)$ to $(\R  \ltimes \heis_{2n+1}/ \mathfrak{a}^+, q_0)$, with
$q_0=\begin{pmatrix}
0&-1&0\\
-1&0&0\\
0&0&I_n
\end{pmatrix}$. 
Let $\varphi \in \Aut(\R \ltimes \Heis_{2n+1})$ be the Lie group isomorphism that corresponds to $\phi$: it acts by an homothety on $\Heis_{2n+1}$. Then consider the map $\psi: (K \times \R)\ltimes \Heis_{2n+1} \to (K \times \R)\ltimes \Heis_{2n+1}$, $\psi(k,g)=(k, \varphi(g))$. Since the action of $K$ on $\Heis_{2n+1}$ commutes with the homotheties of $\Heis_{2n+1}$, $\psi$ is a Lie group automorphism. Furthermore, the induced map $\psi: (G_\rho /I, g) \to (G_\rho /I, g_0)$ is an isometry, where $g$ (resp. $g_0$) is the $G_\rho$-invariant metric on $G_\rho /I$ defined by $q$ (resp. $q_0$). 
\end{proof}
\begin{remark}\label{Remark: q depends only on aL}
One can check in the proof of Proposition \ref{Proposition: when G preserves a metric on G/I} that replacing $L$ by $\alpha L$, for $\alpha \neq 0$, gives a metric isometric to that determined by $L$. So the metric does not depend on the choice of the derivation $L$ generating the $1$-parameter group $\rho(t)=e^{tL}$. 
\end{remark}
\begin{corollary}[\textbf{Statement without compact factor}]\label{Corollary: homog plane wave as $1$-parameter group, without K}
Let $R_\rho=\R \ltimes_\rho \Heis_{2n+1}$ with $\rho(t)=e^{t L}$, $L \in \Der(\heis_{2n+1})$, and set $I=A^+$. Then  $R_\rho/I$ admits a $R_\rho$-invariant Lorentzian metric if and only if the matrix $D$ appearing in $L$ is definite.  
\end{corollary}
\begin{proof}
The `only if' part is proved in Proposition \ref{Proposition: G preserves metric on G/I implies D definite}. The `if' part follows from Proposition \ref{Proposition: when G preserves a metric on G/I} and Reduction 2 just before.     
\end{proof}
\begin{observation}
In Section \ref{Section: Synthetic description of plane waves}, we established a correspondence between plane wave metrics on $X:=\R^{n+1} \times I$, $I \subset \R$ an open interval, whose $\F$-foliation is given by $\R^{n+1} \times \{u\}$, $u \in I$, and some curves in $\Aut(\Heis_{2n+1})$. It appears from this section that the homogeneous ones correspond exactly to those curves which are $1$-parameter groups in $\Aut(\Heis_{2n+1})$!
\end{observation}

\noindent \textbf{Reduction 2.} 
By Remark \ref{Remark: q depends only on aL}, we can suppose that $L(z)= \delta z$, with $\delta \in \{0,1\}$. The reduction in Lemma \ref{Lemme: Blau Derivation form} below will be useful in Section \ref{Section: Global Brinkmann coordinates}, when proving the existence of global Brinkmann coordinates. 
\begin{lemma}\label{Lemme: Blau Derivation form}
Let $G_\rho$ with $\rho(t)=e^{t L}$ as in (\ref{Eq: New derivation L}). Up to isomorphism of $G_\rho$, the derivation $L$ takes the form
$L=\begin{pmatrix}
\delta I+F&B&0\\ I&F&0\\0&0&\delta\end{pmatrix}$
with respect to the decomposition  $\heis_{2n+1}=\a^+\oplus \a^- \oplus \z$, where $F$ is antisymmetric, $B$ is symmetric and $\delta \in \{0,1\}$. 
\end{lemma}
\begin{proof}
The matrix $\begin{pmatrix}
   I_n & P \\
   0 & I_n
\end{pmatrix}$, with $P:=\delta I_n-\frac{A+A^\top}{2}$,
is in $\mathsf{Sp}_{2n}(\R)$. Let $\psi$ be the corresponding (Lie algebra) automorphism of $\heis_{2n+1}$. One can check that for $L$ as in (\ref{Eq: New derivation L}), $\psi L \psi^{-1}$ has the given form. 
\end{proof}
From the study carried out in this section, we deduce the following theorem  characterizing completely the identity component of the isometry group of a $1$-connected non-flat homogeneous plane wave. 
\begin{theorem}\label{THEOREM: Isometry group homog case}
Let $(\mathbf{X},g,V)$ be a $1$-connected non-flat homogeneous plane wave. Then
\begin{enumerate}
\item \textbf{Isometry group.} The connected component of $\Isom(\mathbf{X})$ (which is of finite index) is isomorphic to $G_\rho=(\R \times K) \ltimes \Heis_{2n+1}$,  with the $\R$-action given by $\rho(t)=e^{t L}$, $L \in \Der(\heis_{2n+1})$, such that 
     \begin{itemize}
        \item[(a)] \
          $L=\begin{pmatrix}
           F & B & 0\\
           I & F & 0\\
           0 & 0 & 0
        \end{pmatrix}$,  or  
        \item[(b)] \ 
          $L=\begin{pmatrix}
            I+F & B & 0\\
            I & F & 0\\
            0 & 0 & 1
           \end{pmatrix}$
     \end{itemize}
   with respect to the decomposition  $\heis_{2n+1}=\a^+\oplus \a^- \oplus \z$, where $F$ is antisymmetric and $B$ is symmetric. Moreover, $K$ is the connected (closed) subgroup of $\O(n)$ with Lie algebra $\mathfrak{k}=\{E \in \mathfrak{o}(n), [E,F]=[E,B]=0 \}$, acting on $\Heis_{2n+1}$ trivially on the center and by the standard action on $A^+$ and $A^-$. \medskip
 
 Finally, $\mathbf{X}$ identifies with $X_\rho=G_\rho/K\ltimes A^+$. Moreover, the codimension $1$ foliation $\mathcal{F}$ tangent to $V^\perp$ is given by the (left) action of $K \ltimes \Heis_{2n+1}$, and is preserved by the  action of $G_\rho$.\medskip

\item \textbf{Homogeneous plane waves vs $1$-parameter groups in $\Aut(\Heis_{2n+1})$.} Let $G_\rho$ as in item (1), and consider the homogeneous space $X_\rho=G_\rho/I$, with $I=K \ltimes A^+$. Then $X_\rho$ admits a unique (up to isometry of $X_\rho$) $G_\rho$-invariant Lorentzian metric, which is necessarily a plane wave metric (which may be flat).    \medskip
\item \textbf{Completeness.} If the $\R$-action is as in case (a), the plane wave is geodesically complete. Otherwise, i.e. case (b), they are incomplete. 
\end{enumerate}
\end{theorem}
\begin{proof}
Items (1) and (2) are a summary of the previous results of Section \ref{Section: Isometry group}. Item (3) follows from Remark \ref{Remark: completeness from Brinkman coordinates}. 
\end{proof}

\begin{remark}\label{Remark: B=0}
In the notations of Theorem \ref{THEOREM: Isometry group homog case}, Cahen--Wallach spaces are characterized by the fact that $B$ is non-degenerate and $FB=BF$. On the other hand, $X_\rho$ can be flat, and the flat case occurs exactly when $B=0$.
\end{remark}

\begin{example}
		According to Remark \ref{Remark: B=0}, it enough to choose $BF=FB$ with $B$ non-degenerate to get a Cahen-Wallach space, i.e. an indecomposable symmetric plane wave. Setting $B=I_2$ and $F$ to be a rotation   of angle $\frac{\pi}{2}$,  Theorem \ref{THEOREM: Isometry group homog case} implies that $K=\SO(2)$ in this case. This gives a clear evidence that the more symmetric the plane wave is, the larger $K$ becomes.
\end{example}

\section{Global Brinkmann coordinates}\label{Section: Global Brinkmann coordinates}
In this section, we show the existence of global Brinkmann coordinates for a $1$-connected non-flat homogeneous plane wave. We distinguish two types of spaces, given in the next definition.
\begin{definition}
Let $\mathbf{X}$ be a $1$-connected homogeneous non-flat plane wave.\\
It will be said to be of $(a)$-type if the action of the derivation $L$ is of type $(a)$. They are geodesically complete.\\
It will be said to be of $(b)$-type if the action of the derivation $L$ is of type $(b)$. They are geodesically incomplete. 
\end{definition}
Let $(\mathbf{X},g,V)$ be a $1$-connected homogeneous plane wave. Recall that $\Heis_{2n+1} \subset \Isom(\mathbf{X})$ acts transitively on the $\F$-leaves. Therefore $\mathbf{X}$ is isometric to $X=S/I$, where $S= \R \ltimes \Heis_{2n+1}$, $I=A$  an abelian subgroup of $\Heis_{2n+1}$ intersecting the center only in the identity, and $\rho(t) =e^{tL}$ with $L \in \Der(\heis_{2n+1})$.
Denote by $\a\subset\heis_{2n+1}$ the Lie algebra of $A$. Replacing $L$ by $L+\ad(h)$ for a suitable element $h\in \heis_{2n+1}$, we may assume that $\heis_{2n+1}=\R\oplus \R^{2n}$ as a vector space, where $L(\R^{2n})=\R^{2n}$ and $\a\subset \R^{2n}$. Then $\heis_{2n+1}=\heis(\omega)$ and $\Heis_{2n+1}=\Heis(\omega)$ for a symplectic form $\omega$ on $\R^{2n}$. We write elements of the Heisenberg group as $(z,\xi)$, $z\in\R$, $\xi\in\R^{2n}$. Elements of the $\R$-factor of $S=\R\ltimes\Heis_{2n+1}$ are denoted by $u$. Let $\a'$ be a vector space complement of $\mathfrak{a}$ in $\R^{2n}$. We have the following,
\begin{proposition}
The map 
\begin{eqnarray*}
\phi:\quad \quad X&\longrightarrow & \R\times\mathfrak{a}'\times\R\\
(u,z,\xi)\cdot A &\longmapsto & (v:=z+\textstyle \frac12\omega(x,y),y,u),
\end{eqnarray*}
where $\xi=x+y$ for $x\in \mathfrak{a}$ and $y\in\mathfrak{a}'$, is a diffeomorphism.
\end{proposition}
\begin{proof} The map is well defined. Indeed, let $\hat x$ be in $A$. Then 
$(u,z,\xi)\cdot(0,0,\hat x)=(u, z+\textstyle\frac12\omega(\xi,\hat x),\xi+\hat x)$, which is mapped to $(z+\textstyle\frac12\omega(\xi,\hat x)+\frac12\omega(x+\hat x,y),y,u)$. Now we use that $\xi=x+y$ and that $\omega(x,\hat x)=0$ since $\a$ is an abelian subalgebra. Thus we get the same image as for $(u, z,\xi)$. Now, we show that $\phi$ is a bijective local diffeomorphism. For injectivity, if we assume that $\phi((u,z,\xi)A)=\phi(u'\cdot (z',\xi')A)$ a straight forward computation we get $y=y', u=u'$ and $z'=z-\textstyle\frac12\omega(x-x',y)$, so we get that $(u',z',\xi')=(u,z,\xi)(0,0,x'-x)\in (u, z,\xi)A$. The surjectivity is clear since, we can cover the $y$ and $u$ coordinates. For the $v$ coordinates it is enough to choose $\xi$ such that $\omega(x,y)=0$ so $z$ can cover the rest. To see that $\phi$ is a local diffeomorphism we compute the differential at an arbitrary point $p:=(u,z,\xi)A$, where we can assume that $\xi \in \a'$. We have that $$\gamma_t:= (t,0,0) \cdot (u,z,\xi) = ((t+u,z, \exp(tL)\xi),$$ and let $p_1$, $p_2$ be the projection on $\a$, $\a'$ respectively. We compute the image of $\partial u$ by $d_p\phi$. Namely,
\begin{eqnarray*}
d_p\phi(\partial_u) = && \odv{\phi(\gamma_t)}{t}_{t=0} \\ 
= && \odv{\left(t+u,z+\frac{1}{2}\omega\left(p_1(\exp(tL)\xi), p_2(\exp(tL)\xi)\right), p_2(\exp(tL)\xi)\right)}{t}_{t=0}  \\ 
= && \left(1, \frac{1}{2}\omega\left(p_1(L\xi), p_2(\xi) \right) + \omega\left(p_1(\xi), p_2(L\xi)\right), p_2(L\xi)\right) \\ 
= && \left(1, \frac{1}{2}\omega\left(p_1(L\xi), \xi) \right), p_2(L\xi)\right).
\end{eqnarray*}
Hence, we see that $d_p\phi(\partial_u)=\partial_u + \cdots$. We do an analogues computation with $\partial_\xi$ and $\partial_z$ we deduce that the differential has a full rank, the proposition follows.  
\end{proof} 

Now suppose that $x^1,\dots,x^n$ are coordinates of $\a\cong\R^n$, $y^1,\dots,y^n$ are coordinates of $\a'$ and that $\a'$ is chosen such that the symplectic form on $\a\oplus\a'\cong\R^{2n}$ is given by $\omega((x,y), (x',y'))= (x,y)J {x'\choose y'}$ for $J=\left(\begin{matrix}
0&I\\-I&0\end{matrix}\right)$.
We identify the tangent space of $X$ at $ o:=eA\in S/A$ with $\mathfrak{s}/\a$. The metric $g$ on $X$ induces an $\ad_{\a}$-invariant scalar product $\langle\cdot, \cdot\rangle$ on $\mathfrak{s}/\a$. Recall that $\a'=\{0\}\times\R^n$. Thus $\mathfrak{s}/\a\cong \R\oplus\a'\oplus\R=\R\oplus\R^n\oplus\R$ with coordinates $z,y,u$. We determine $\langle\cdot, \cdot\rangle$ with respect to these coordinates. The $\ad_{\a}$-invariance implies that up to isometries of $X$, $\langle\cdot, \cdot\rangle$ equals 
$$\begin{pmatrix}
0&0&1\\0&I_n&0\\1&0&0\end{pmatrix}.$$

\begin{theorem} \label{Pa} If $\mathbf{X}$ is a $1$-connected homogeneous plane wave of $(a)$-type, then $\mathbf{X}$ is isometric to $\R\times\R^n\times \R$ with coordinates $(v,\bar y, u)$ and metric given by
$$2du dv+d\bar y^2+\big(e^{u F} B e^{- uF}\big)_{ij}\bar y^i \bar y^j  {d  u^2}.$$
Conversely, any such metric on $\R\times\R^n\times \R$ defines a $1$-connected homogeneous plane wave of $(a)$-type. 
\end{theorem}
\begin{proof}
Let $\mathbf{X}$ be of $(a)$-type . Then a direct calculation shows 
$$(\phi^{-1})^* g= 2du dv +2(Fy \boldsymbol{\cdot} dy )du +dy^2 + (Fy \boldsymbol{\cdot} Fy+ y \boldsymbol{\cdot}By)du^2,$$
where `$\boldsymbol{\cdot}$' denotes the standard euclidean product. The coordinate transformation 
$$\bar y:= e^{uF}y$$ gives 
\[dy^j=-(F e^{-uF})^j_k \bar y ^k du + (e^{-uF})^j_k d\bar y ^k\ =\ -(Fy)^j du +(e^{-uF})^j_k d\bar y ^k.\]
This implies 
\begin{eqnarray*}
(Fy \boldsymbol{\cdot} dy )du&=& \textstyle \sum_j (Fy)^j dy^jdu\\
&=&\textstyle \sum_j (Fy)^j (-(Fy)^j du +(e^{-uF})^j_k d\bar y ^k)du\\
&=& \textstyle -(Fy \boldsymbol{\cdot} Fy) du^2 +\sum_j (Fy)^j(e^{-uF})^j_k d\bar y ^k du\\[0.5ex]
(dy)^2&=&\textstyle \sum_j (dy^j)^2 \ =\ \sum_j (-(Fy)^j du +(e^{-uF})^j_k d\bar y ^k)^2\\
&=& \textstyle (Fy\boldsymbol{\cdot}Fy) du^2 - 2\sum_j (Fy)^j(e^{-uF})^j_k d\bar y^k du +\sum_j (e^{-uF})^j_k (e^{-uF})^j_l d\bar y ^k d\bar y ^l\\
&=& \textstyle (Fy\boldsymbol{\cdot}Fy) du^2 - 2\sum_j (Fy)^j(e^{-uF})^j_k d\bar y^k du +d\bar y^2
\end{eqnarray*}
since $e^{-uF}$ is an orthogonal matrix because $F$ is antisymmetric. Furthermore,
\[(y \boldsymbol{\cdot}By)du^2= (e^{-uF}\bar y)\boldsymbol{\cdot} (Be^{-uF}\bar y)\,  du^2\\
= \big( \bar y \boldsymbol{\cdot} (e^{uF}Be^{-uF}\bar y) \big)\, {d u^2}\]
since $(e^{-uF})^\top = e^{uF}$. Consequently,
\begin{eqnarray*}
(\phi^{-1})^*g &=& 2du dv-2 \textstyle (Fy \boldsymbol{\cdot} Fy) du^2 +2\sum_j (Fy)^j(e^{-uF})^j_k d\bar y ^k du\\
&&+\textstyle (Fy\boldsymbol{\cdot}Fy) du^2 - 2\sum_j (Fy)^j(e^{-uF})^j_k d\bar y^k du +d\bar y^2\\
&&+ (Fy \boldsymbol{\cdot} Fy)du^2 +\big( \bar y \boldsymbol{\cdot} (e^{uF}Be^{-uF}\bar y) \big)\, {d u^2}\\
&=& 2d u d v +d\bar y^2+\big( \bar y \boldsymbol{\cdot} (e^{uF}Be^{-uF}\bar y) \big)\, {d u^2}  
\end{eqnarray*}
hence the existence of the given coordinates. The converse is straightforward and follows from the fact that there is a Heisenberg action transitive on each $u$-level, and a $1$-parameter flow of isometries acting transversally to the $u$-levels, given by: for any $s \in \R, s \cdot (v,\Bar{y},u) = (v,e^{sF} \Bar{y}, u+s)$.
\end{proof}

\begin{theorem}\label{pb} If $\mathbf{X}$ is a $1$-connected homogeneous plane wave of $(b)$-type , then $\mathbf{X}$ is isometric to $\R\times\R^n\times \R_{>0}$ with coordinates $(\bar v,\bar y,\bar u)$ and metric given by 
$$2d\bar u d\bar v+d\bar y^2+\big(e^{\ln(\bar u)F} B e^{-\ln(\bar u)F}\big)_{ij}\bar y^i \bar y^j \frac {d \bar u^2}{\bar u^2}.$$
Conversely, any such metric on $\R\times\R^n\times \R$ defines a $1$-connected homogeneous plane wave of $(b)$-type.
\end{theorem}
\begin{proof} Let $\mathbf{X}$ be of $(b)$-type . As in the proof of Theorem ~\ref{Pa}, a direct calculation shows 
$$(\phi^{-1})^* g= 2du dv +2(Fy \boldsymbol{\cdot} dy )du +dy^2 + (2v+Fy \boldsymbol{\cdot} Fy+ y \boldsymbol{\cdot}By)du^2.$$
 The coordinate transformation $\tau$ 
$$\bar u :=e^{-u}, \quad \bar y:= e^{uF}y, \quad
\bar v:= -e^{u}v$$ gives 
\begin{eqnarray*}
dv&=&e^{-u} \bar v du-e^{-u} d\bar v \ =\ -v du  -e^{-u} d\bar v,\\
dy^j&=& -(Fy)^j du +(e^{-uF})^j_k d\bar y ^k,\\
du&=& -e^u d\bar u\,.
\end{eqnarray*}
This implies 
\[ du dv = du(-v du  -e^{-u} d\bar v)\ =\ -v du^2 +d\bar u d\bar v.\]
Moreover, as in the proof of Theorem ~\ref{Pa},
\begin{eqnarray*}
(Fy \boldsymbol{\cdot} dy )du&=&  \textstyle -(Fy \boldsymbol{\cdot} Fy) du^2 +\sum_j (Fy)^j(e^{-uF})^j_k d\bar y ^k du\\[0.5ex]
(dy)^2&=&\textstyle (Fy\boldsymbol{\cdot}Fy) du^2 - 2\sum_j (Fy)^j(e^{-uF})^j_k d\bar y^k du +d\bar y^2
\end{eqnarray*}
since $e^{-uF}$ is an orthogonal matrix because $F$ is antisymmetric. Furthermore,
\begin{eqnarray*}(y \boldsymbol{\cdot}By)du^2&= &(e^{-uF}\bar y)\boldsymbol{\cdot} (Be^{-uF}\bar y)\,  e^{2u}d\bar u^2\\
&= &\big( \bar y \boldsymbol{\cdot} (e^{uF}Be^{-uF}\bar y) \big)\,  \frac{d\bar u^2}{\bar u ^{2}}.
\end{eqnarray*}
We obtain
\begin{eqnarray*}
(\tau^{-1})^*(\phi^{-1})^*g &=& 2(-v du^2 +d\bar u d\bar v)-2 \textstyle (Fy \boldsymbol{\cdot} Fy) du^2 +2\sum_j (Fy)^j(e^{-uF})^j_k d\bar y ^k du\\
&&+\textstyle (Fy\boldsymbol{\cdot}Fy) du^2 - 2\sum_j (Fy)^j(e^{-uF})^j_k d\bar y^k du +d\bar y^2\\
&&+ (2v+Fy \boldsymbol{\cdot} Fy)du^2 +\big( \bar y \boldsymbol{\cdot}(e^{uF}Be^{-uF}\bar y) \big)\,  \frac{d\bar u^2}{\bar u ^{2}}\\
&=& 2d\bar u d\bar v +d\bar y^2+\big(\bar y \boldsymbol{\cdot}(e^{-\ln (\bar u)F}Be^{\ln(\bar u)F}\bar y) \big)\,  \frac{d\bar u^2}{\bar u ^{2}}
\end{eqnarray*}
hence the existence of the given coordinates. The converse again follows from the fact that there is a Heisenberg action transitive on each $u$-level, and a $1$-parameter flow of isometries acting transversally to the $u$-levels, given by: for any $s \in \R, s \cdot (v,\Bar{y},u) = (s v,e^{\log(s^{-1}) F} \Bar{y}, s^{-1}u)$.
\end{proof}

\begin{remark}[The flat case]
		As stated in Remark \ref{Remark: isometry group of flat PW}, a $1$-connected homogeneous flat plane wave is isometric to either the whole Minkowski space or a half-Minkowski space bounded by a lightlike plane. Both spaces have global Brinkmann coordinates.
\end{remark}

\begin{remark}\label{Remark: completeness from Brinkman coordinates}
    Item $(3)$ in Theorem \ref{THEOREM: Isometry group homog case} can be  deduced from the existence of global Brinkmann coordinates. Indeed, with these coordinates in hand, the completeness of $(a)$-type metrics follows from \cite[Proposition 3.5]{sanchez-candela}. And the incompleteness of $(b)$-type metrics follows from the fact that the geodesic parameter of a geodesic transversal to $\F$ coincides with the $u$-coordinate (up to affine change).
\end{remark}
\section{$C^2$-inextendibility}\label{Section: C^2 maximality}
We say that a pseudo-Riemannian manifold $(X,g)$ is $C^k$-inextendible (or maximal) if there is no isometric  $C^k$-embedding of $(X,g)$ as a proper open subset of some pseudo-Riemannian manifold $(Y,h)$.
\begin{theorem}\label{Theorem: Maximality}
    Let $\mathbf{X}$ be a $1$-connected non-flat homogeneous plane wave. Then $\mathbf{X}$ is $C^2$-inextendible.
\end{theorem}
We begin with the following lemma, whose Lemma \ref{Lemma: rho(Heis)=Heis} is a special case :
\begin{lemma}\label{Lemma: rho(Heis)=Heis in L}
    Let $\varphi$ be a locally faithful representation of the Heisenberg group $H:=\Heis_{2n+1}$ in $\mathsf{L}(n)=(\R \times \O(n))\ltimes \Heis_{2n+1}$. Then $\varphi(H)=\Heis_{2n+1} \subset \mathsf{L}(n)$. In particular, $\varphi$ is faithful.
\end{lemma}
\begin{proof}
Recall that the $\R$-action on $\Heis_{2n+1}$ in $\mathsf{L}(n)$ is given by the derivation $$\begin{pmatrix}I+F & 0 & 0\\0 & F & 0\\0&0&1\end{pmatrix},$$
written in the basis $(\a^+, \a^-, \z)$ of $\heis_{2n+1}$.
    Let $H_0:=\varphi(H)$. Let $p$ be the projection from $\mathsf{L}(n)$ to the $\R$-factor. If $p(H_0)$ is trivial, we conclude using Lemma \ref{Lemma: rho(Heis)=Heis}. If not, we could write $H_0=\R\ltimes H_1$, where $H_1\subset \O(n)\ltimes \Heis_{2n+1}$ (observe that $H_1$ is connected). 
    From the proof of Lemma \ref{Lemma: rho(Heis)=Heis}, we see that either $H_1$ is contained in $\Heis_{2n+1}$ or $\dim (H_2:= H_1\cap \Heis)=2n-1$. 
    If $n>1$, then $2n-1\geq n+1$, hence $\z\subset H_2$. But  $[L,\z]=\z$, and this implies that $H_0$ is not nilpotent. In the case $n=1$, $H_0\cong \R\ltimes \R^2$, where the $\R$-action is nilpotent and semisimple, hence trivial, but the algebraic multiplicity of the $0$-eigenvalue of $L$ is $1$. 
\end{proof}
\begin{proof}[Proof of Theorem \ref{Theorem: Maximality}]
We can suppose that $\mathbf{X}$ is simply connected and of type (b) (for otherwise, it is complete). Suppose $\mathbf{X}$ is not maximal, and let $\phi : \mathbf{X} \to Y$ be a $C^2$ embedding into a bigger Lorentzian manifold $Y$. For simplicity of notation, we identify $\phi(\mathbf{X})$ with $\mathbf{X}$.
\subsection*{Step 1: extending the codimension 1 foliation}
Codimension 1 geodesic ``laminations''  of pseudo-Riemannian manifolds, and more generally manifolds with an affine  connection,   satisfy a general  Lipschitz regularity principle. Let $U$ be a convex open subset of a manifold $Y$ endowed with a connection $\nabla$, that is geodesic segments  exit uniquely between points of $U$ (and are contained in $U$). Let $Z$ be a subset of $U$ endowed with a partition $\mathcal F$, such that for any $x \in Z$,
 $\mathcal F (x)$ is a   geodesic hypersurface, closed in $U$. In the case $\dim Y = 2$, $\mathcal F(x)$ is a geodesic, closed in $U$, and so has its endpoints in $\partial U$.  The same picture applies in higher dimension. Then, $\mathcal F$ extends to the closure
 of $Z$ in $U$ and is locally Lipschitz (see \cite[Theorem 7.1]{zeghib2004lipschitz} for details).

In our case, let us apply this fact to  $U$, a convex neighborhood of a point $p_0 \in \overline{\mathbf{X}} \smallsetminus \mathbf{X}$. We consider the foliation on $U \cap \mathbf{X}$ given by the connected components of the intersection of the $\mathcal F$-leaves (defined in $\mathbf{X}$) with $U$. Actually, we  forget the (global) $\mathcal F$ foliation of $\mathbf{X}$ and replace it by this foliation of $\mathbf{X} \cap U$, and keep the same notation $\mathcal F$ for it.  So it extends as a foliation of $\overline{\mathbf{X}} \cap U$. 
 
Let $F_0$ be the leaf of $p_0$, and  $\gamma:  [0, \epsilon[ \to Y$, $\gamma(0) = p_0$,  a geodesic transversal to $F_0$. Then, it is transversal to all $\mathcal F$-leaves near $p_0$. Let $ \tau = \gamma([0, \epsilon[)$, 
then $\tau^\prime:= \tau  \cap U$  is open in $\tau \cong [0, \epsilon[$. 
Assume  $\tau^\prime$ contains an interval $\gamma(]0, \beta[)$, that is a neighborhood of $p_0$ in $\tau$ except $p_0$, contained in $\mathbf{X}$. 
Then the saturation of $\gamma([0, \beta[)$  by $\mathcal F$ will be a codimension 0 submanifold with boundary $F_0$. 
 
Let us show that we can always choose $p_0$ in order to have this property  (of a neighbourhood in $\tau$ except its endpoints, which is contained in $\mathbf{X}$). 
For this, consider a connected component of $\tau^\prime$.  It has the form $\gamma (]\alpha, \beta[)$, $\alpha, \beta \in ]0, \epsilon[$.  It is clear that we if  we  replace $p_0$ by $p_0^\prime= \gamma (\alpha)$, and the geodesic $\gamma$ by its restriction  to $]\alpha, \beta[$, then $\gamma (]\alpha, \beta[) \subset \mathbf{X}$, and so near $p_0^\prime$,  $\overline{\mathbf{X}} \cap U $ is a regular open set with boundary the $\mathcal F$-leaf of $p_0^\prime$.

\subsection*{Step 2: extending small isometries of $G_{\rho}$} Consider the regular open set $\overline{\mathbf{X}} \cap U$ with boundary the $\mathcal F$-leaf of $p_0$, which we denote by $F_0$. The extended foliation will be denoted $\overline{\F}$. 
Here, we prove that the infinitesimal action of $\mathfrak{g}_{\rho}$ on $\overline{\mathbf{X}} \cap U$ extends to $F_0$ by fixing $F_0$. Let $f\in G_{\rho}$ in a small neighborhood of the identity. Let $\gamma_0$ be a geodesic though $p_0$ transversal to $\overline{\F}$. Let
$p \in \gamma_0 \cap U$. For $p$ close enough to $F_0$, there is an open neighborhood $B$ of $\gamma_0'$ at $p$ in $T_p \mathbf{X}$ such that all geodesics from $B$ intersect $F_0$. Hence they form an open neighborhood of $p_0$ in $F_0$.  
Let $X \in \mathfrak{g}_{\rho}$. For $f^t=\exp(t X)$ close enough to the identity, $f^t(\gamma)$ such that $\gamma'(0) \in B$ intersects $F_0$, and we have $f^t(\gamma \cap F_0) = f^t(\gamma) \cap F_0$. This defines a (unique) smooth extension of $f^t$ to a neighborhood of $p_0$. This extension preserves the affine lightlike geometry of $F_0$. Hence we have a representation $\psi: \mathfrak{g}_\rho=(\R \oplus \mathfrak{k})\ltimes \heis_{2n+1} \to \mathfrak{l}=(\R \oplus \mathfrak{o}(n)) \ltimes \heis_{2n+1}$, where $\mathfrak{l}$ denotes the Lie algebra of the affine lightlike group.  This representation is faithful; indeed, an isometry acting trivially on a hypersurface (even degenerate) is trivial. Moreover, we have $\psi(\mathfrak{z}) \subset \mathfrak{z}$, i.e. $\psi$ sends the center of $\heis_{2n+1} \subset \mathfrak{g}_\rho$ to the center of $\heis_{2n+1} \subset \mathfrak{l}$ (indeed, by the above extension process, the lightlike parallel vector field extends to $F_0$ to a lightlike parallel vector field, necessarily generated by the action of the center of $\heis_{2n+1}$).

\subsection*{End of the proof}
Let $\psi: \mathfrak{g}_\rho=(\R \oplus \mathfrak{k})\ltimes \heis_{2n+1} \to \mathfrak{l}=(\R \oplus \mathfrak{o}(n)) \ltimes \heis_{2n+1}$ be a faithful representation. We have by Lemma \ref{Lemma: rho(Heis)=Heis in L} that $\psi(\heis_{2n+1}) = \heis_{2n+1}$. In our special situation, we can prove this also using the fact that $\psi(\mathfrak{z}) \subset \mathfrak{z}$.
Indeed, let $L$ (resp. $L'$) be a generator of the $\R$-factor in $\mathfrak{g}_\rho$ (resp. in $\mathfrak{l}$). Let $z$ (resp. $z'$) be a generator of the center of $\heis_{2n+1} \subset \mathfrak{g}_\rho$ (resp. of $\heis_{2n+1} \subset \mathfrak{l}$). The adjoint action of any element of $\mathfrak{k}\ltimes \heis_{2n+1} \subset \mathfrak{g}_\rho$ (resp. of $\mathfrak{o}(n) \ltimes \heis_{2n+1} \subset \mathfrak{l}$) on $z$ (resp. on $z'$) is trivial, whereas the adjoint action of $L'$ on $z'$ is non-trivial. So $\psi(\mathfrak{k}\ltimes \heis_{2n+1}) \subset \mathfrak{o}(n) \ltimes \heis_{2n+1}$, hence $\psi(\heis_{2n+1}) = \heis_{2n+1}$ by Lemma \ref{Lemma: rho(Heis)=Heis}. 

Let $\mathfrak{i}=\mathfrak{k} \ltimes \a^+$.  So $\mathfrak{g}_\rho/\mathfrak{i}$, with its $\ad_\mathfrak{i}$-invariant Lorentzian scalar product (recall that the latter is unique, up to isometry), is isometric to $\phi(\mathfrak{g}_\rho)/\psi(\mathfrak{i})$, with its $\ad_{\psi(\mathfrak{i})}$-invariant scalar product.  We have 
\begin{align*}
    L=\begin{pmatrix}
			I+F_0 & B\\
			I & F_0
		\end{pmatrix} \mathrm{\;\;and\;\;}
   L'=\begin{pmatrix}
			I & 0\\
			0 & 0
		\end{pmatrix}, 
\end{align*}
in the decomposition $\heis_{2n+1} = \a^+ \oplus \a^- \oplus \mathfrak{z}$,  with $F_0 \in \mathfrak{o}(n)$. Moreover $\psi(L)=L'+F_1= \begin{pmatrix}
			I+F_1 & 0\\
			0 & F_1
		\end{pmatrix}$, with $F_1 \in \mathfrak{o}(n)$. 
 We claim that $\psi(\a^+) \cap \a^- = {0}$. Before proving the claim, let us see that it   leads to a contradiction. Indeed, consider a new decomposition $\heis_{2n+1}=\psi(\a^+)\oplus \a^- \oplus \mathfrak{z}$, where $\psi(\a^+)$ is a Lagrangian of $\heis_{2n+1}$. Then the matrix of $\psi(L)$ in this decomposition is  $\begin{pmatrix}
			I+F_1 & 0\\
			* & F_1
		\end{pmatrix},$ which is conjugate in $\Aut(\heis_{2n+1})$ to $\begin{pmatrix}
			I+F_1' & 0\\
			I & F_1'
		\end{pmatrix},$ with $F_1' \in \mathfrak{o}(n)$. But then this means that the space is flat (see Remark \ref{Remark: B=0}), which contradicts our assumption on $\mathbf{X}$. Let us now prove the claim. Suppose that there is $h \in \a^+$ such that $\psi(h) \in \a^-$.  On the one hand, $\psi([L,h])=\psi(L(h))=\psi(h+F_0(h)+h')$, where $h'\in \a^-$ satisfies $[h,h']=\alpha z$, $\alpha \neq 0$.  On the other hand, using that $\psi$ is a Lie algebra isomorphism and that $\psi(h) \in \a^-$, we get $\psi([L,h])=[\psi(L),\psi(h)]=F_1(\psi(h)) \in \a^-$.  It follows that $\psi(h+F_0(h)+h') \in \a^-$, hence $\psi(F_0(h)+h') \in \a^-$. As a consequence, $\psi([h,F_0(h)+h'])=[\psi(h),\psi(F_0(h)+h')]=0$. 
  But we also have $[h,F_0(h)+h']=\alpha z$, with $\alpha \neq 0$, hence $\psi([h,F_0(h)+h'])= \alpha \psi(z) \neq 0$. This ends the proof. 
\end{proof}
\appendix
\section*{Appendix. Action of the Heisenberg group in Brinkmann coordinates}\label{app}

Consider $\mathbf{X}=\R^{n+1} \times I$ with coordinates $(v,x=(x_1,\dots,x_n),u)$, equipped with the metric
	\begin{equation} \label{Eq2: Brinkmann coordinates}	
	g = 2 dv du + x {}^\top S(u) x\,  du^2 + \Sigma_{i=1}^n dx_i^2,
	\end{equation}  
where $S(u)$ is a symmetric matrix. This is the metric of a plane wave in Brinkmann coordinates, where the lightlike parallel vector field is given by $V:=\partial_v$. In what follows, we compute the isometries of $(\mathbf{X},g)$ acting trivially on the $u$-coordinate. Denote this group by $G_0$.\\
Let $S: u \in I \mapsto S(u)$ be the curve of symmetric matrices defining the metric $g$.
\begin{fact}\label{Appendix-Fact: computation of G_0 in Brinkmann coordinates}
An isometry $\varphi \in G_0$ has the following form
\begin{align*}
		\varphi  
		\begin{pmatrix}
			v \\
			x\\
			u
		\end{pmatrix}&=
		\begin{pmatrix}
			v - \langle\alpha'(u), Ax + \frac{1}{2} \alpha(u)\rangle + c  \\
			Ax+\alpha(u)  \\
			u
		\end{pmatrix},
\end{align*} 
where $A \in C_{\O(n)}(S)$, and $\alpha$ is a solution of the differential equation
\begin{align*}
    \alpha''(u)-S(u)\alpha(u)=0, \;\forall u \in I.
\end{align*}
\end{fact}
In the following, we prove the above fact.
We know that $V=\partial_v$ is a parallel lightlike vector field, and the foliation tangent to $V^\perp$ is given by the levels of the $u$-coordinate. Let $\varphi$ be an isometry of $g$ acting trivially on the $u$-coordinate. Set $\varphi=(\varphi_v, \varphi_1,\dots,\varphi_n, u)$. 
First, since $V$ is the only lightlike vector tangent to the $u$-levels, $\varphi$ preserves the $V$-direction, i.e. $\varphi_*(V)=a V$ for some $\a \in \R^*$. This gives $\frac{\partial \varphi_i}{\partial v}=0$ for any $i$, hence the $\varphi_i$ components do not depend on $v$.

Now, at each $u_0$-level, $\varphi$ induces a map from $V^\perp / V = \R^n$ to itself given by $\Bar{\varphi}_{u_0}=(\varphi_1,\dots,\varphi_n)_{\vert \{u=u_0\}}$, which preserves the induced Euclidean scalar product $\sum dx_i^2$. It follows that $\Bar{\varphi}_{u_0}$ is a Euclidean transformation (that depends only on $u_0$), hence the existence for any $u \in I$ of $A(u) \in \O(n)$ and $\alpha(u) \in \R^n$ such that $\Bar{\varphi}(v,x,u) = A(u)x + \alpha(u)$. \\
\\
\textbf{Step 1:} $A(u)$ is a constant curve in $\O(n)$. 
\begin{proof}[Proof of Step 1]
The equality  $g(\varphi_* \partial_u, \varphi_* \partial_i)=g(\partial_u, \partial_i)=0$ reads 
\begin{align}\label{Eq0}
\frac{\partial \varphi_v}{\partial x_i} + \langle A(u)(e_i), A'(u)x+\alpha'(u)\rangle =0.    
\end{align}
Then, writing $\frac{\partial}{\partial x_k}\frac{\partial \varphi_v}{\partial x_i} = \frac{\partial}{\partial x_i}\frac{\partial \varphi_v}{\partial x_k}$ yields $\langle A(u)(e_i), A'(u)(e_k) \rangle = \langle A(u)(e_k), A'(u)(e_i) \rangle$ for any $i,k \in \{1,\dots,n\}$, proving that the matrix $B(u):=A(u)^\top A'(u)$ is symmetric. But since $A(u)$ is a curve in $\O(n)$ we also have that $B(u)$ is antisymmetric, implying $B(u)=0$ for any $u$, hence the result.  
\end{proof}
\noindent In what follows, $A(u)$ will be then simply denoted by $A$.\\
\\
\textbf{Step 2:} From Step 1 we have $$\Bar{\varphi}(v,x,u)=Ax+ \alpha(u)$$ for some $A \in \O(n)$ and some function $\alpha$ depending only on $u$. 

Writing $g(\varphi_* \partial_u, \varphi_* \partial_v)=g(\partial_{u}, \partial_v)=1 $ gives $a \frac{\partial \phi_u}{\partial u}=1$, hence $a=1$, and then $\frac{\partial \varphi_v}{\partial v}=1$. On the other hand, it follows from (\ref{Eq0}) that $\frac{\partial \phi_v}{\partial x_i}=- \langle A(e_i),\alpha'(u) \rangle$. 
From this one obtains 
\begin{align}\label{varphi v}
\varphi_v=v - \langle \alpha'(u), Ax \rangle + l(u),    
\end{align}
with $l$ some function depending only on $u$.\\
Finally, using $g(\varphi_* \partial_u, \varphi_* \partial_u)=g(\partial_u, \partial_u)$ together with (\ref{varphi v}) and the fact that $S(u)$ is symmetric, gives for any $x \in \R^n$ and $u \in I$
\begin{align*}
x^\top B(u)x - 2 x^\top A^\top Y(u) + L(u)=0,
\end{align*}
with 
\begin{align*}
B(u):&=A^\top S(u) A - S(u),\\
Y(u):&=\alpha''(u)-S(u) \alpha(u),\\    
L(u):&=\langle \alpha(u), S(u)\alpha(u) \rangle + \langle \alpha'(u), \alpha'(u) \rangle + 2 l'(u).
\end{align*}
From this, and using that $A$ is invertible, we obtain that for any $u \in I$
\begin{align*}
    &S(u)=A^\top S(u) A,\\
    &\alpha''(u)-S(u) \alpha(u)=0,\\
    &l'(u)=-\frac{1}{2} (\langle \alpha', \alpha'\rangle - \langle \alpha, \alpha''\rangle).
\end{align*}
As a consequence, $A$ commutes with $S(u)$ for all $u \in I$, $\alpha$ is a solution of the linear differential equation 
\begin{align}\label{diff equ alpha''(u)-S(u) alpha(u)=0}
    \alpha''(u)-S(u) \alpha(u)=0,
\end{align}
and $l(u)=- \frac{1}{2}\langle \alpha', \alpha \rangle + c$, with $c \in \R$. Finally, we get
\begin{align}\label{varphi form when A=I}
		\varphi
		\begin{pmatrix}
			v \\
			x\\
			u
		\end{pmatrix}&=
		\begin{pmatrix}
			v - \langle\alpha'(u), Ax + \frac{1}{2} \alpha(u)\rangle + c  \\
			Ax+\alpha(u)  \\
			u
		\end{pmatrix},
\end{align}
where $A \in C_{\O(n)}(S)$, $\alpha$ is a solution of the differential equation (\ref{diff equ alpha''(u)-S(u) alpha(u)=0}) and $c \in \R$. Conversely, one can check that any application of this form is an isometry of $g$.
\bigskip

\bibliographystyle{abbrv}
\bibliography{Bibliography}

\begin{thebibliography}{10}

\bibitem{Besse}
A.~L. Besse.
\newblock {\em Einstein manifolds}.
\newblock Classics in Mathematics. Springer-Verlag, Berlin, 2008.
\newblock Reprint of the 1987 edition.

\bibitem{Blau}
M.~Blau and M.~O'Loughlin.
\newblock Homogeneous plane waves.
\newblock {\em Nuclear Physics B}, 654(1):135--176, 2003.

\bibitem{blumenthal}
R.~A. Blumenthal.
\newblock Transversely homogeneous foliations.
\newblock {\em Ann. Inst. Fourier (Grenoble)}, 29(4):vii, 143--158, 1979.

\bibitem{cahen1970lorentzian}
M.~Cahen and N.~Wallach.
\newblock Lorentzian symmetric spaces.
\newblock {\em Bulletin of the American Mathematical Society}, 76(3):585--591,
  1970.

\bibitem{sanchez-candela}
A.~M. Candela, J.~Flores, and M.~Sanchez.
\newblock On general plane fronted waves. geodesics.
\newblock {\em General Relativity and Gravitation}, 35(4):631--649, 2003.

\bibitem{choquetbruhat}
Y.~Choquet-Bruhat and R.~Geroch.
\newblock Global aspects of the {C}auchy problem in general relativity.
\newblock {\em Communications in Mathematical Physics}, 14:329--335, 1969.

\bibitem{duncan1989homogeneous}
D.~C. Duncan and E.~C. Ihrig.
\newblock Homogeneous spacetimes of zero curvature.
\newblock {\em Proceedings of the American Mathematical Society},
  107(3):785--795, 1989.

\bibitem{ehrlich1992gravitational}
P.~E. Ehrlich and G.~G. Emch.
\newblock Gravitational waves and causality.
\newblock {\em Reviews in Mathematical Physics}, 4(02):163--221, 1992.

\bibitem{Leis}
W.~Globke and T.~Leistner.
\newblock Locally homogeneous pp-waves.
\newblock {\em J. Geom. Phys.}, 108:83--101, 2016.

\bibitem{Gro}
M.~Gromov.
\newblock Rigid transformations groups.
\newblock In {\em G\'{e}om\'{e}trie diff\'{e}rentielle ({P}aris, 1986)},
  volume~33 of {\em Travaux en Cours}, pages 65--139. Hermann, Paris, 1988.

\bibitem{Article2}
M.~Hanounah, I.~Kath, L.~Mehidi, and A.~Zeghib.
\newblock Topology and dynamics of compact plane waves.
\newblock {\em Journal f{\"u}r die reine und angewandte Mathematik (Crelle's
  Journal)}, 280:87--113, 2025.

\bibitem{Leis-Schlie}
T.~Leistner and D.~Schliebner.
\newblock Completeness of compact {L}orentzian manifolds with abelian holonomy.
\newblock {\em Math. Ann.}, 364(3-4):1469--1503, 2016.

\bibitem{Ross}
D.~Marolf and S.~F. Ross.
\newblock Plane waves: {T}o infinity and beyond!
\newblock {\em Classical and Quantum Gravity}, 19(24):6289, 2002.

\bibitem{nomizu}
K.~Nomizu.
\newblock On local and global existence of {K}illing vector fields.
\newblock {\em Annals of Mathematics}, pages 105--120, 1960.

\bibitem{thurston2022geometry}
W.~P. Thurston.
\newblock {\em The Geometry and Topology of Three-Manifolds: With a Preface by
  Steven P. Kerckhoff}, volume~27.
\newblock American Mathematical Society, 2022.

\bibitem{zeghib2004lipschitz}
A.~Zeghib.
\newblock Lipschitz regularity in some geometric problems.
\newblock {\em Geometriae Dedicata}, 107:57--83, 2004.

\end{thebibliography}
	
\end{document}